\documentclass{amsart}
\usepackage{graphicx}
\newtheorem{thm}{Theorem}[section]

\newtheorem{lem}[thm]{Lemma}

\theoremstyle{definition}
\newtheorem{defn}[thm]{Definition}
\theoremstyle{remark}
\newtheorem{rem}[thm]{Remark}
\numberwithin{equation}{section}

\numberwithin{equation}{section}

\newcommand{\To}{\longrightarrow}

\begin{document}

\title[The Existence of Quasimeromorphic Mappings in Dimension 3]{The Existence of Quasimeromorphic Mappings in Dimension 3}%
\author{Emil Saucan}%
\address{Technion and Ort Braude}%
\email{semil@tx.technion.ac.il}%

\thanks{This paper represents part of the authors Ph.D. Thesis written under the supervision of Prof. Uri Srebro}%

\date{20.9.2003.}
\begin{abstract}
We prove that a Kleinian group $G$ acting upon $\mathbb{H}^{3}$ admits a non-constant $G$-automorphic function,
even if it has torsion elements, provided that the orders of the elliptic (i.e torsion) elements are uniformly
bounded. This is accomplished by developing a technique for meshing distinct fat triangulations while preserving
fatness. We further show how to adapt the proof to higher dimensions.
\end{abstract}
\maketitle
\section{Introduction}
The object of this article 
is the study of the existence of $G$-automorphic quasimeromorphic mappings (in the sense of Martio and Srebro --
see [MS1]) $f:\mathbb{H}^n \rightarrow \widehat{\mathbb{R}^n}$, $\widehat{\mathbb{R}^n} = \mathbb{R}^n \bigcup
\,\{\infty\}$; i.e. such that
 \begin{equation} f(g(x))= f(x) \quad ;
  \qquad \forall x \in \mathbb{H}^n ; \forall g \in G \; ;
 \end{equation}
were $G$ is Kleinian group acting upon $\mathbb{H}^n$.
\\ Our principal goal is to prove the following:

\begin{thm} Let $G$ be a Kleinian group with torsion acting upon $\mathbb{H}^n, \, n \geq 3$. \\ If the elliptic elements (i.e. torsion elements) of $G$ have
uniformly bounded orders,
\\then there exists a non constant $G$-automorphic quasimeromorphic mapping
\\$f: \mathbb{H}^n \rightarrow \widehat{\mathbb{R}^n}.$
\end{thm}

In this paper we restrict ourselves to the proof of the theorem in the classical case (i.e. $n = 3$) only. This
restriction is motivated by two reasons: (a) the proof in the $3$-dimensional case employs mainly elementary tools
and (b) it develops and uses a technique for for meshing distinct fat triangulations while preserving fatness,
technique that is relevant in Computational Geometry and Mathematical Biology. The proof of the general case is
presented in \cite{s1} and it is based upon a more general result concerning the existence of fat triangulations
for manifolds with boundary -- see \cite{s2}.

The question whether quasimeromorphic mappings exist was originally posed by Martio and Srebro in \cite{ms1} ;
subsequently in \cite{ms2} they proved the existence of fore-mentioned mappings in the case of co-finite groups
i.e. groups such that $ Vol_{hyp} (H^n/ G) < \infty  $ (the important case of geometrically finite groups being
thus included). Also, it was later proved by Tukia  (\cite{tu}) that the existence of non-constant
quasimeromorphic mappings (or qm-maps, in short) is assured in the case when G acts torsionless upon
$\mathbb{H}^n$. Moreover, since for torsionless Kleinian groups $G$, $\mathbb{H}^n/G$ is a (analytic) manifold,
the next natural question to ask is whether there exist qm-maps $f:M^n \rightarrow \widehat {\mathbb{R}^n} $;
where $M^n$ is an orientable $n-$manifold. The affirmative answer
to this question is due to K.
 Peltonen (see \cite{pe}); to
be more precise she proved the existence of qm-maps in the case when $M^n$ is a connected, orientable $
\mathcal{C^\infty}$\!-Riemannian manifold.
 \newline In contrast with the above results it was proved by Srebro (\cite{sr}) that, for any $n\geq 3$, there exists a Kleinian group $G \;\rhd \! \!\! \!<
 \mathbb{H}^n$ s.t. there exists no non-constant, $G-$automorphic function  $f:\mathbb{H}^n\rightarrow \mathbb{R}^n$.
 More precisely, if $G$ (as above) contains elliptics of unbounded orders (with non-degenerate fixed set), then $G$ admits no non-constant
$G-$automorphic qm-mappings.
\newline Since all the existence results were obtained in constructive manner by using the classical "Alexander
trick" (See \cite{al}), it is only natural that we try to attack the problem using the same method. For this
reason we present here in succinct manner Alexander's method: One starts by constructing a suitable triangulation
(Euclidian or hyperbolic) of $\mathbb{H}^n$ or of $M^n = \mathbb{H}^n/G$. Since $\mathbb{H}^n$ and $M^n$ are
orientable, an orientation consistent with the given triangulation  (i.e. such that two given $n$-simplices having
a $(n-1)$-dimensional face in common will have opposite orientations) can be chosen. Then one quasiconformally
maps the simplices of the triangulation into $\widehat{\mathbb{R}^n}$ in a chess-table manner: the positively
oriented ones onto the interior of the standard simplex in $\mathbb{R}^n$ and the negatively oriented ones onto
its exterior. If the dilatations of the $qc$-maps constructed above are uniformly bounded (as is the case of
compact manifolds $M^n$), then the resulting map will be quasimeromorphic.
\newline The dilatations of each of the $qc$-maps above is dictated by the proportions of the respective simplex (see \cite{tu} , \cite{ms2}), and since the dilatation is to
be uniformly bounded, we are naturally  directing our efforts in the construction of a {\bf fat} triangulation,
where:
\begin{defn} A k-simplex $\tau \subset \mathbb{R}^n$ (or $\mathbb{H}^n$); $2 \leq k \leq (n-1)$
is $f${\it -fat} if there exists $f \geq 0$ s.t. the  ratio $\frac{r}{R} \geq f$; where $r$ denotes the radius of
the inscribed sphere of $\tau$ ({\sf inradius}) and $R$ denotes the radius of the circumscribed sphere of $\tau$
({\sf circumradius}).
\\ A triangulation (of a submanifold of $\mathbb{R}^n$ or $\mathbb{H}^n$) $\mathcal{T} = \{ \sigma_i \}_{i\in \bf I}$ is $f${\it -fat} if all its simplices are f-fat.
\\ A triangulation $\mathcal{T} = \{ \sigma_i \}_{i\in \bf I }$ is {\it fat} if there exists $f \geq 0$ s.t. all its simplices are $f${\it-fat}; $\forall i \in \bf{I}. $
\end{defn}
The idea of the proof of Theorem 1.1. is first to build two fat triangulations: $\mathcal{T}_{e}$ of a certain
closed neighbourhood $\overline{N}_e$ of the fixed set of $G$ in $\mathbb{H}^{n}$\,; and $\mathcal{T}_c$ of
$\mathbb{H}^{n}\: \backslash \;\overline{N_e}$; and then to "mash" the two triangulations into a new
triangulation, while retaining their fatness\footnote{\,But see also \cite{ca1}, \cite{ca2}.}.
\\The first triangulation is constructive and is based upon the geometry of the elliptic transformations. The existence of the second triangulation is assured by Peltonen's result.
Unfortunately, these two triangulations are not $G$-invariant, so they are unsuited for our purpose of building a
$G$-automorphic function. However, they induce fat triangulations:
 ${\mathcal{T}}_{e}^\ast$ on $(\overline{N_e} \cap \mathbb{H}^{n})/G$\,, and $\overline{\mathcal{T}_{p}}$ on
$\overline{M_{c}} = (\mathbb{H}^{n}\setminus N_e)/\,G$, where $M_{c}$ is a differential manifold with boundary
$\partial M_{c} = (\partial \overline{N_e} \bigcap \mathbb{H}^{n})/\,G$.
\\Fortunately, we are provided with a ready made method of mashing triangulations,
so we can direct our efforts towards the task of "fattening" the simplices of the "intermediate zone"; task which
will be carried out in Section 4.
\\ The "mashing" method mentioned above is based on a result and, even more, on the technique used in its proof, due to Munkres:

\begin{thm} Let $M^{n}$ be a $\mathcal{C}^{r}$-manifold with boundary. Then any $\mathcal{C}^{r}$-triangulation of $\partial M^{n}$ can be extended to a
$\mathcal{C}^{r}$-triangulation of $M^{n}$, $1\leq r\leq \infty$.
\end{thm}
Because of its importance to our own construction, we shall present the basic idea of the proof of Theorem 1.3. in
the next Section.
\\ This paper is organized as follows: in Section 2 we present the necessarily background on elliptic
transformations and present in a nutshell the main techniques we employ: the Alexander trick, Peltonen's method
and the Proof of Munkres' Theorem. In Section 3 we show how to choose and triangulate the closed neighbourhood of
the $\overline{N}$ of the fixed set of $G$, and how to select the "intermediate zone" where the two different
triangulations overlap. Section 4 is dedicated to the main task of fattening the common triangulation. Finally, in
Section 5 we indicate the way of adapting our construction to higher dimensions.

\section{Preliminaries}

\subsection{Elliptic Transformations} We shall restrict ourselves mainly to the $3$-dimen-sional case, for, as we
have already stated, this will be the direction in which our main efforts will be directed.
\\ Let us first recall the basic definitions and notations: A transformation $f\in Isom(H^n)$, $f\neq Id$ is called {\it elliptic} if $(\exists)\,m\geq 2 $ s.t. $f^{m} = Id$,
and $m$ is called the {\it order} of $f$. In the $3$-dimensional case the {\it fixed point set} of f, i.e. $
Fix(f) = \{ x\in H^{n} | f(x) = x \}$, is a hyperbolic line and will be denoted by $A(f)$ - the {\it axis of f}.
If $A$ is an axes of an elliptic of order $m$, then $A$ is called an $m-axes$.
\\  If the discrete group $G$ is acting upon $\mathbb{H}^3$, then by the discreteness of $G$, there exists no accumulation point
of the elliptic axes in $\mathbb{H}^3$. Moreover, if $G$ contains no elliptics with intersecting axes, then the
distances between the axes are, in general, bounded from bellow. To be more precise, the following holds:

\begin{thm}[\,\cite{gm1}] Let $G$ be a discrete group $G$ acting upon $\mathbb{H}^3$, and let $f,g \in G$ be s.t.
$ord(f) \geq 3$ or $ord(g) \geq 3$;  and s.t. $A(f) \cap A(g) = \emptyset$.
\\ Then $\exists\, \delta >0$, that is independent of $G,f,g$ s.t.
  \begin{equation}
     dist_{hyp}(A(f),A(g)) \geq \delta \,;
  \end{equation}
where $dist_{hyp}$ denotes the hyperbolic distance in $\mathbb{H}^3$.
\end{thm}

 It is extremely important to notice that the results above do not include the case when all the elliptic transformations of $G$ are of order $2$, since they
depend intrinsically upon:
\begin{thm} [\cite{j} , \cite{ah}]
 $G < H^{n}$ is discrete iff \ $<f,g>$ is discrete; $\forall f,g \in G$.
\end{thm}
The theorem above is of little avail in the case of two {\it half-turns} (that is elliptics of order $2$); since
any two half-turns (in $H^{3}$) generate a discrete group. (See \cite{teza} for a proof of this  "folkloric"
result from an unpublished paper by J\o rgensen.) Indeed, examples of discrete groups of isometries of hyperbolic
$3$-space can be constructed, such that the distances between the axes of the order $2$ elliptics are not bound
from bellow (see \cite{teza}).

In the presence of {\it node points} (i.e. intersections of axes) the situation is more complicated. Fortunately,
there are only a few types of such possible intersections - for the orders of the elliptic axes meeting at a node
point must satisfy certain conditions (determined by the Euler number for the orbifold\footnote{\, See \cite{th1},
\cite{th2} for the definition and the necessary proofs.}); namely, the possible local situations in orbifold  are:
either (a) {\em Dihedral}, i.e. of the type $(2,2,n)\,, n\geq 2$; or one of the following exceptional types: (b)
{\em Tetrahedral}, i.e. of type $(2,2,3)$\,; \;(c) {\em Octahedral}\,, i.e. of type $(2,3,4)$\,;\; (d) {\em
Icosahedral}, i.e.  of type $(2,3,5)$;\; (see Fig. 1).

 \begin{figure}[h]
 \begin{center}
\includegraphics[scale=0.25]{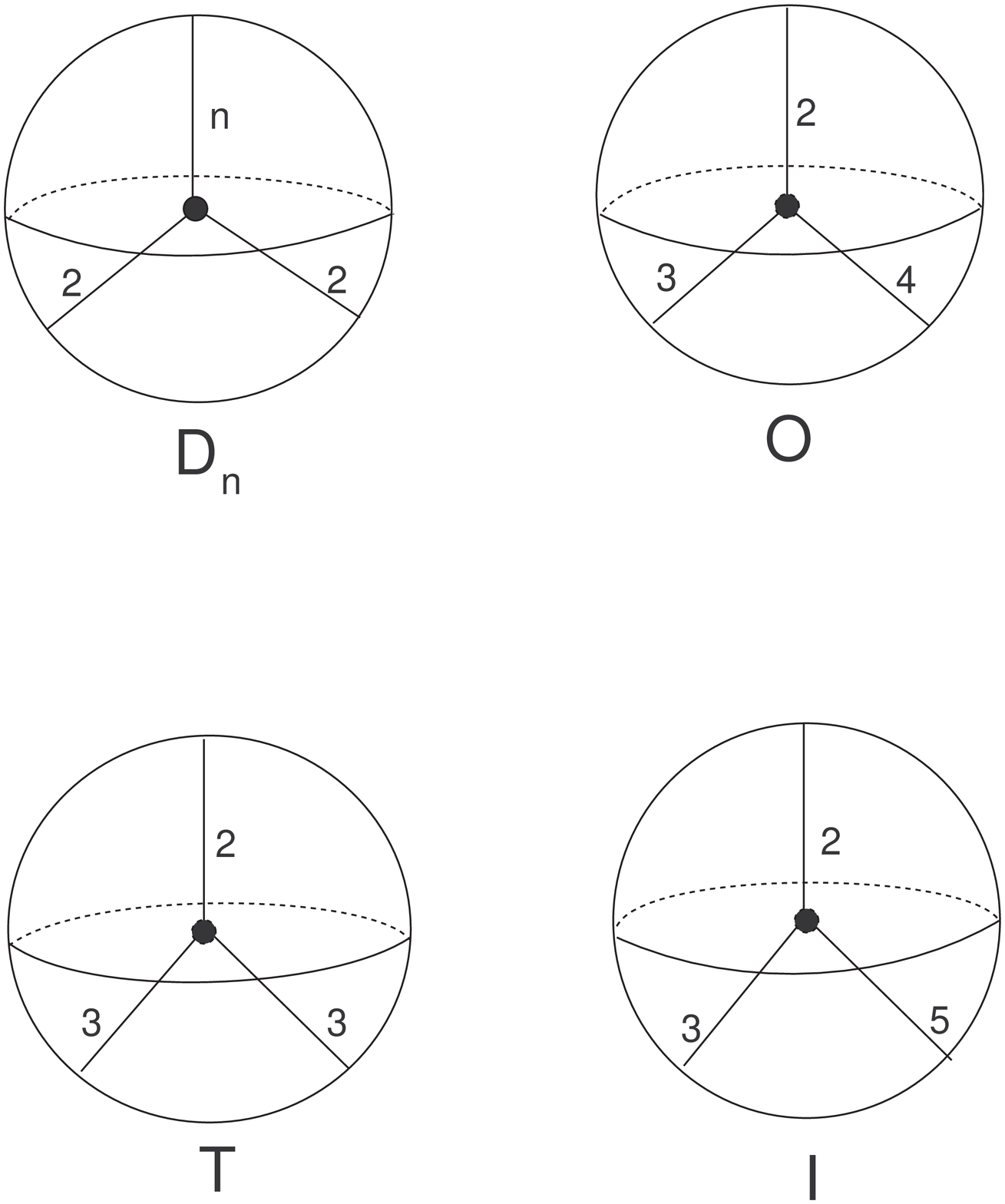}
\end{center}
 \caption{}
\end{figure}

\begin{rem}
The reduced number of possibilities is rather fortunate, for the computation of the distances between node points
is more difficult than that of distances between disjoint axes. (For more specific information in this direction
of study see \cite{dm} and \cite{med}\footnote{\, and also \cite{teza}.}; and more recently  \cite{gm1} ,
\cite{gm2} and \cite{gmmr}).
\end{rem}

Since our main interest lies in Kleinian groups acting upon $\mathbb{H}^3$ whose elliptic elements have orders bounded from above, the following theorem is highly relevant:

\begin{thm}[\cite{fm}]
Let $G$ be finitelly generated Kleinian group acting on $H^3$\,.
Then the number of conjugacy classes of elliptic elements is finite.
\end{thm}
For a sketch of an alternative proof of this Theorem see Appendix.

\begin{rem}
The Theorem above is not true for groups of isometries of $H^n\,, n\geq 4$\,; indeed there exist counterexamples,
one due to Mess and Feighn (\cite{fm}) and another due to Kapovitch and Potyagailo (\cite{kp}). It should be
remarked that both of the examples cited above produce (albeit different) conjugacy classes of elliptics of the
{\em same} order. Considering this and the goal of our investigation, the following recent result is highly
relevant:
\begin{thm}[\cite{h}]
    There exists a discontinuous group $\Gamma < Isom(H^{n})\,, \:n \geq 4$\,; such that:
    \\(i)\;$\Gamma$ contains elliptics of arbitrary large orders
    \\ and
    \\(ii) $Vol(N_{\varepsilon}(M_{\Gamma})) < \infty$\footnote{\,That $\Gamma$ is {\em almost} geometrically finite.}
    , where $M_{\Gamma} = H_{\Gamma}/\Gamma$, and
    $H_{\Gamma}$ denotes the convex core (in $H^{n}$) of the limit set $\Lambda(\Gamma) \subset
    H^{n}$, and $N_{\varepsilon}$ represents the $\varepsilon$-neigbourhood of $M_{\Gamma}$.
\end{thm}
\end{rem}

\subsection{Alexander's Trick}
The technical ingredient in Alexander's trick is the following Lemma (which we formulate for $\mathbb{R}^{3}$
only, but which readily generalizes to higher dimensions):
\begin{lem} (\cite{ms1}, \cite{pe})
Let $\mathcal{T}$ be a fat triangulation of $M \subset \mathbb{R}^3$, and let $\tau,\sigma \in
\nolinebreak[4]\mathcal{T},\; \tau = (p_1,p_2,p_3,p_4), \,\sigma = (q_1,q_2,q_3,q_4)$; and denote $|\tau| = \tau
\cup int\,\tau$.
\\ Then there exists a sense-preserving homeomorphism $h = h_{\tau}: |\tau| \rightarrow \widehat{\mathbb{R}^{3}}$
s.t.
\begin{enumerate}
\item $h(|\tau|) = |\sigma|$, \,if\; $det(p_1,p_2,p_3,p_4) > 0$
\\ and
\\ $h(|\tau|) = \widehat{\mathbb{R}^{3}} \setminus| \sigma|$, \,if\; $det(p_1,p_2,p_3,p_4) < 0$.
\item $h(p_i) = q_i, \; i=1,\ldots,4.$
\item $h|_{\partial|\sigma|}$ is a $PL$-homeomorphism.
\item $h|_{int|\sigma|}$ is quasiconformal.
\end{enumerate}
\end{lem}

\begin{rem}
The branching set of $h$ is the $1$-skeleton of the triangulation.
\end{rem}

\subsection{Peltonen's Technique} Peltonen's method is an extension of one due to Cairns, developed in order to triangulate $\mathcal{C}^{2}$-compact manifolds
(\cite{ca3}). It is based on the subdivision of the given manifold into a closed cell complex generated by a
Dirichlet (Voronoy) type partition whose vertices are the points of a maximal set that satisfy a certain density
condition. We give below a sketch of the Peltonen's method, refereing the interested reader to the authoritative
\cite{pe} for the full details.\footnote{\, A rather detailed exposition of the main steps of the proof one can be
find in \cite{teza}}
\\The construction devised by Peltonen consists of two parts:
\\{\em Part 1} This part of the proof proceeds in two steps:
\\  \hspace*{0.3cm}{\em Step A} We build an exhaustation $\{E_i\}$ of $M^n$, generated by the pair $(U_i,\eta_i)$, where:
\begin{enumerate}
\item $U_i$ is th relatively compact set $E_i \setminus E_{i-1}$ and
\item $\eta_i$ is a number that controls the "fatness" of the simplices of the triangulation of $E_i$, that will be constructed in Part 2, such that  they don't differ to much
on adjacent simplices, i.e.:
\\ (i) The sequence $(\eta_i)_{i\geq1}$ descends to $0$\,;
\\ (ii) $2\eta_i \geq \eta_{i-1} \,.$
\end{enumerate}
 \hspace*{0.3cm}{\em Step B}
\begin{enumerate}
\item Produce a maximal set $A$, $|A| \leq \aleph_0$, s.t. $A \cap U_i$ satisfies:
\\ (i) a density condition, and
\\ (ii) a "gluing" condition (for $U_i, U_{i=1}$).
\item Prove that the Dirichlet complex $\{\bar{\gamma}_i\}$ defined by the sets $A_i$ is a cell complex and
every cell has a finite number of faces (so it can be triangulated in a standard manner).
\end{enumerate}
{\em Part 2} Consider first the dual complex $\Gamma$ and prove that it is a Euclidian simplicial complex with a
"good" density, then   project $\Gamma$ on $M^n$ (using the normal map). Finally, prove that the resulting complex
can be triangulated by fat simplices.

\subsection{Munkres' Theorem} The basics steps in the proof\footnote{\, A much more detailed exposition of the proof is given in \cite{teza}. For the full proof one should consult, of course, the original study of Munkres \cite{mun}.  } of Theorem 1.3. are as follows:
\\ a) Prove that you can triangulate a smooth manifold without boundary in the following way: approximate $M^{n}$
locally by a locally finite Euclidian triangulation, by means of the secant map (see \cite{mun}, p. 90). Modify
these local triangulations coordinate chart by chart, so they will be $PL$-compatible wherever they overlap. To
extend the triangulation globally, we work in $\mathbb{R}^{n}$, by using the coordinate charts and maps. Here
again we have to approximate the given triangulation by a  $PL$-map, s.t. the given triangulation and the one we
produce will be compatible.
\\ b) Triangulate a product neighbourhood  $\mathcal{P}(\partial M^{n})$ of  $\partial M^{n}$, $\mathcal{P}(\partial M^{n}) \subset
M^{n}$, in a standard way, and mash it together with the triangulation of the non-bounded manifold $int\,  M$, by
using the same method as above.
\\ It is important to emphasize that for most of the process the technique
sketched above not only preserves the fatness of the simplices, but actually takes care that the said fatness will
occur.

\section{Constructing and Intersecting Triangulations}

\subsection{Geometric Neigbourhoods}
If there are no elliptics with intersecting axes, a standard choice for a regular neighbourhood of an $m$-axes
will be -- for obvious geometric reasons -- a doubly-infinite regular hyperbolic $m$-prism (henceforth called a
{\it geometric} neighbourhood), and a fundamental domain will be a prismatic "slice", i.e. a fundamental region
for the action of $C_m$ on $\{m\} \times A(f)$, where $\{m\}$ denotes the regular hyperbolic polygon with $m$
sides and $C_m$ denotes the cyclic group of order $m$, i.e. the rotation group of $\{m\}$.  In order to
triangulate the geometric neighbourhood, we divide it in a finite number of radial strata of equal width $\varrho
= \delta/\kappa_0$, and further partition it into "slabs" of equal hight $h$. Each prismatic fundamental region
thus obtained naturally decomposes into three congruent tetrahedra, generating a $C_m$-invariant triangulation of
the geometric neigbourhood. (See Fig. 2 for a representation in the ball model of $H^3$ in the case $m = 4$.) The
fatness of the triangulation of the geometric neighbourhood thus depends upon $\delta$, $\kappa_0$ and $h$,
enabling one to control the initial fatness of the geometric triangulation by means of the parameters $\varrho$
and $h$.

\begin{figure}[h]
\begin{center}
\includegraphics[scale=0.25]{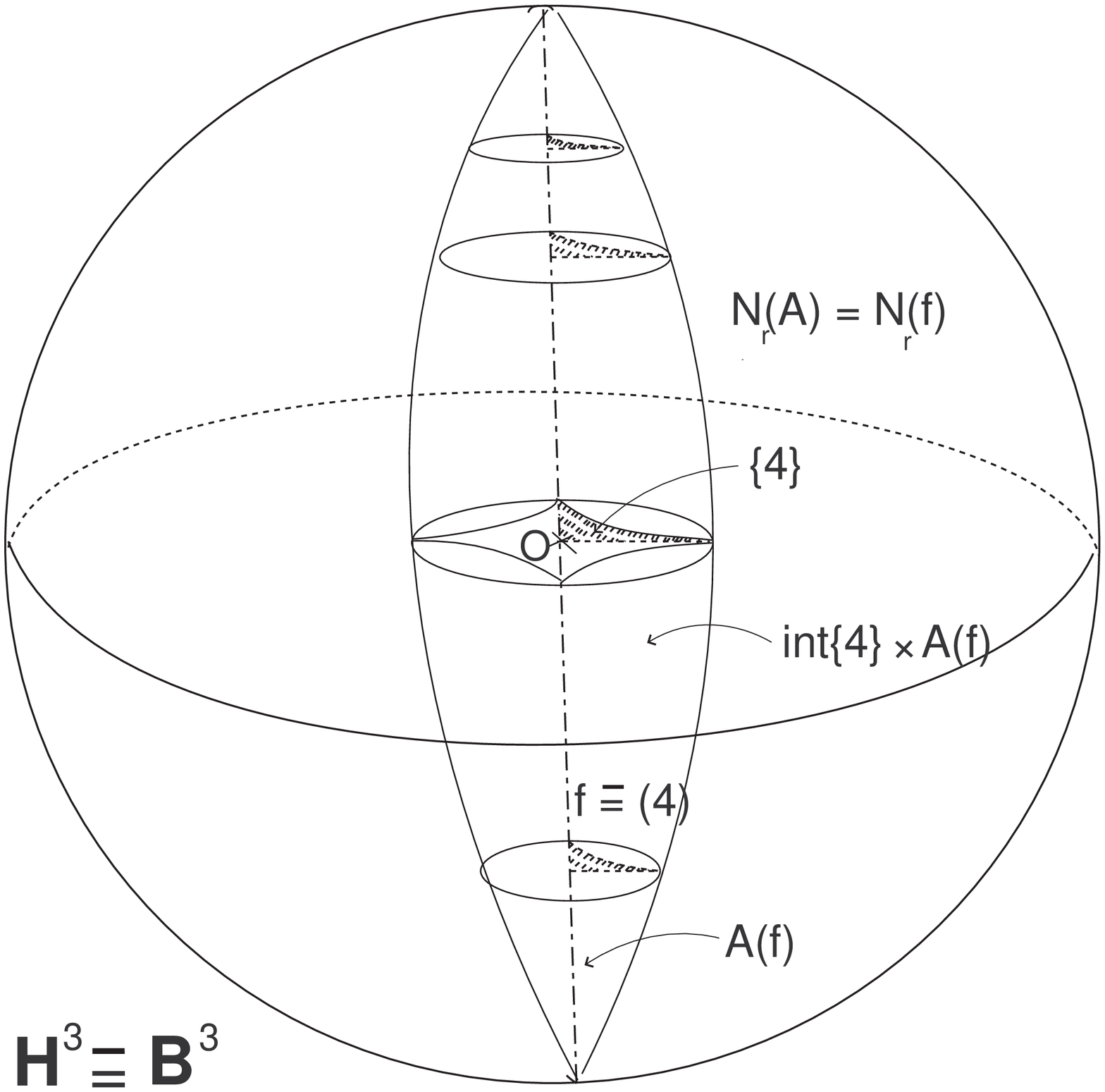}
\end{center}
\caption{}
\end{figure}

\begin{rem}
One can easily modify the construction above and produce  instead of a $C_m$-invariant triangulation, a
$D_m$-invariant one, where $D_m$ denotes the {\em dihedral} group of order $m$, i.e. the full-symmetry group of
$\{m\}$, thus allowing one to consider groups that contain orientation-reversing isometries of $\mathbb{H}^3$.
\end{rem}

For the choice of geometric neigbourhoods for the node points, the natural choice is that of  an Archimedean solid
which is a natural carrier of the symmetry group of the desired type: $D_n\,, T,\,O$ or $I$ (or rather for its
spherical counterpart -- see \cite{cox})\nolinebreak[4].

\subsection{Mashing Triangulations} We start by showing first how to construct the desired fat triangulation and how to
produced the quasimeromorphic mapping ensuing from it in the basic case of groups who's elliptic elements axes do
not intersect. Moreover, let us presume here that there exist a least an elliptic element of order $\geq 3$. Since
$G$ is a discrete group, $G$ is countable so we can write $G = \{g_j\}_{j \geq 1}$ and let $\{f_i\}_{i \geq 1}
\subset G$ denote the set of elliptic elements of $G$. The steps in building the fat desired fat triangulation are
as follows:

\begin{enumerate}

\item Let $N_i = \{x \in \mathbb{H}^3 \,|\, dist_{hyp}(A_i,x) < \delta/4\}$, where $\delta/4$ is the constant provided by Theorem 2.1, and where
$A_i = A(g_i)$. We also put: $N_e = \bigcup \raisebox{-0.7em}{\hspace{-0.5cm}\tiny $i\in \mathbb{N}$}\!N_{i}$.

\item Consider the following quotients: $N_i^{\ast} = (\overline{N}_i \cap \mathbb{H}^3)/\,G$, $i \geq 1$; and \newline $N_e^{\ast} = (\overline{N}_e \cap \mathbb{H}^3)/G$.

\item Let $\mathcal{T}_e$ denote a $G$-invariant fat triangulation of $N_e$ that induces a $G$-invariant fat triangulation $\mathcal{T}_i$
of $N_i$, and let $\mathcal{T}_e^{\ast}, \mathcal{T}_i^{\ast}$ denote the fat triangulations induced upon
$N_e^{\ast}, \, N_i^{\ast}$ by $T_e$ and $T_i$, respectively.

\item Denote by $\mathcal{T}_p$ the fat triangulation of $M_c = (\mathbb{H}^3/\,G)\,\backslash \,N_e^{\ast} = (\mathbb{H}^3\, \backslash\, \overline{N}_e)/\,G$ assured by Peltonen's
Theorem. The lift of $\mathcal{T}_p$ to $\mathbb{H}^3$ is a $G$-invariant fat triangulation $\mathcal{T}_c$ of
$\mathbb{H}^3\,\backslash \,N_e$.


\item The desired triangulation $\mathcal{T}^{\ast}$ of $M = \mathbb{H}^3/G$ will consist of  the simplices of $\mathcal{T}_e^{\ast}$, the simplices of
$\mathcal{T}_p^{\ast}$ away from $\partial M_{c} = (\partial \overline{N_e} \bigcap \mathbb{H}^{3})/\,G$, and
other new simplices $\mathcal{T}_b \;\subset\; \partial M_{c} \times [0,1) \; \subset\; M_{c}$\,, that are
constructed by applying the mashing technique of Munkres.

\item Apply Alexander's Trick to $\mathcal{T}^{\ast}$ and get a quasimeromorphic mapping \\$f^{\ast}:\mathbb{H}^3/G \rightarrow
\widehat{\mathbb{R}^3}$.

\item The lifting of $f^{\ast}$ to a mapping $f:\mathbb{H}^3 \rightarrow \widehat{\mathbb{R}^3}$ produces the
desired $G$-invariant quasimeromorphic mapping.

\end{enumerate}

\begin{rem}
The choice of "$\delta/4$" instead of "$\delta$" in the definition of the geometric neighbourhoods $N_i$ is
dictated by the following Lemma:

\begin{lem}(\cite{rat})
      Let $X$ be a metric space, and let $\Gamma < Isom(X)$ be a discontinuous group.
      \\Then, for any $x \in X$ and any $r \in (0,\delta /4)$:
      \[ \pi : B(x,r)/\Gamma_{x} \simeq B(\pi(x),\delta/4) \]
      where: $\Gamma_{x}$ is the stabilizer of $x$, $\pi$ denotes the natural projection,
\\ $\delta := d(x,\Gamma(x) \backslash \{x\})$, and where the metric on $X/\Gamma$ is given by
      \[ d_{\Gamma}([\pi(x)],[\pi(y)]) = d(\Gamma(x),\Gamma(y))\,; \,\forall x,y\in  X/\,\Gamma \,.\]
\end{lem}

\end{rem}

\hspace*{-0.5cm} {\bf Note} \, Instead of the triangulation scheme presented above, scheme that follows closely
the Proof of Theorem 1.3.\,, we could have used in this case the natural triangulation of geometric neigbourhoods
in $H^3$, to devise a simpler method for mashing triangulations, as follows:

\begin{enumerate}
\item Consider again the geometric neighbourhood \[N_i = N_{i,1/4} = \{x \in \mathbb{H}^3 \,|\, dist_{hyp}(A_i,x) < \delta/4\}\] with its natural fat triangulation.
\item Replace the neighbourhoods $N_i = N_{i,1/4}$ by the neighbourhoods \[N_i' = N_{i,3/16} = \{x \in \mathbb{H}^3 \,|\, dist_{hyp}(A_i,x) <
3\delta/16\}\,,\] and triangulate $N_i'$ in such a manner that the simplices of the triangulation of $\partial
N_i'$ are also simplices of the triangulation of $int\,N_i$.
\item Consider instead of $N_e$ the following manifold: $N_e' = \bigcup \raisebox{-0.7em}{\hspace{-0.5cm}\tiny $i\in
\mathbb{N}$}\!N_{i}'$.
\item Replace $M_c$ by $M_c' = (\mathbb{H}^3/\,G)\,\backslash \,N_e'^{\ast} = (\mathbb{H}^3\, \backslash\,
\overline{{N}_e'}\,)/\,G$.
\item The role of the triangulation $\mathcal{T}_p$ is played by $\mathcal{T}_p'$, which consists of the simplices
produced by Peltonen's method and those simplices resulting from those of the geometric triangulation of $T_i =
N_{i,1/4} \backslash  N_{i,3/16}$.
\item The desired triangulation $\mathcal{T'}^{\ast}$ are composed of those of $ N_{i,1/4}$, those of the original
$M_c$ and those obtained by mashing the two triangulations of the tubes $T_i$.
\end{enumerate}

\begin{rem}
The second construction, besides being more simple and geometrically intuitive, reduces more rapidly the original
problem to that of mashing and uniformly fattening two locally finite Euclidian triangulations.
\end{rem}

In the case when there exist intersecting elliptic axes, the following modification of our construction is
required: instead of $\delta$ one has to consider $\delta^* = \min{(\delta, \delta_0)}$, where $\delta_0$
represents the minimal distance between node-points.\footnote{\, See \cite{gm1}, \cite{gm2}, \cite{gmmr}.}

We still have to deal with the case when all the elliptic transformations are half-turns, since, as we have seen,
no minimal distance between the axes can be computed in this case. However, by the discreetness of $G$ it follows
that there is no accumulation point of the axes in $\mathbb{H}^3$. Let $D = \{d_{ij}\,|\, d_{ij} =
dist_{hyp}(A_i,A_j)\}$ denote the set of mutual distances between the axes of the elliptic elements of $G$. Then,
since $G$ is countable, so will be $D$, thus $D = \{d_k\}_{k \geq 1}$. Then  the set of neighbourhoods
$N_e^{\natural} = \bigcup \raisebox{-0.7em}{\hspace{-0.5cm}\tiny $k \in \mathbb{N}$}\!N_{k}^{\natural} = \bigcup
\raisebox{-0.7em}{\hspace{-0.5cm}\tiny $k \in \mathbb{N}$}\!\{x \in \mathbb{H}^3 \,|\, dist_{hyp}(A_k,x) <
\delta/4k\}$ will constitute a proper geometric neigbourhood of $A_G =
\bigcup\raisebox{-0.7em}{\hspace{-0.5cm}\tiny $i \in \mathbb{N}$}\!A_i$. The fatness of the simplices of the
geometric triangulation of $A_G$ can be controlled, as before, by a proper choice of $h$ and $\varrho$.
\begin{rem}
The existence of $N_e^{\natural}$ is easy to justify geometrically if one uses the upper-half space model of
$\mathcal{H}^3$: up to conjugation one can choose $\infty$ to be an accumulation point for $A_G$, therefore the
axes accumulated at this point are represented as parallel Euclidian half-lines, perpendicular to the plane
$\mathbb{R}^2$. Consider a family of disjoint Euclidian cylinders $C_i = \{x\,|\,dist_{Eucl}(A_i,x) < r_i\}$, $C_i
\cap \mathbb{R}^2 = p_i$. Let $l_{i}^{\pm}$ represent two parallel generators of $C_i$, let $h_i$ be the hyperbola
of vertex $p_i$ and asymptotes $l_{i}^{\pm}$, and let $H_i$ the hyperboloid of rotation with axes $A_i$ and
generatrix $l_i^+$. Then $\{int\, H_i\,|\,i \geq 1\}$ represents a proper geometric neighbourhood for the set of
axes accumulating at $\infty$.
\end{rem}

\section{Fattening Triangulations}

\subsection{Preliminaries} We have seen that we reduced the problem to that of "fattening" the intersection of two
$3$-dimensional finite, fat Euclidian triangulations.
\\ We do this piecemeal, first fattening the $2$-simplices, then the $3$-simplices. It is natural to do so, for the following
holds:

\begin{lem}[\cite{pe}]
If an $n$-dimensional simplex is fat, than all its $k$-dimensional faces, $2 \leq k \leq n-1$ are fat.
\end{lem}

In particular, in order that an $n$-dimensional simplex be fat, its $2$-simplices have to be fat. Note that for
triangles, with angles $\alpha\,, \beta\,,\gamma\,,$ and $r$ and $R$ as above, the conditions: $r/R\geq f$ and\,
$\min{\{\alpha\,, \beta\,,\gamma\}} \geq \varphi$, where $\varphi = \varphi(f)$ is an angle depending on $f$, are
equivalent. Thus we start "fattening" $2$-simplices, by ensuring that
\begin{equation}
\min_i{\{\alpha_{i}\,, \beta_{i}\,,\gamma_{i}\}} \geq \frac{\varphi_{0}}{10}\;;
\end{equation}
where the minimum is taken over all the triangles of the resulting triangulation, and where $\varphi_{0}$ is the
minimal angle of the original triangulations -- ensured by their uniform fatness.
\\From now on, let $\mathcal{S} = \{s_{i}\}_{i\in I}$ and $\Sigma = \{\sigma_{j}\}_{j\in J}$ stand for the simplices
of the triangulations of $T_i$ and $M_c'$, respectively.\footnote{\,or: $\partial M_{c} \times [0,1)$ and
$M_{c}$\,.}
\\Since the intersections of tetrahedra can be rather unruly, we simplify the situation
by requiring that the simplicies of one of the triangulations be much smaller than those of the other:

\begin{equation}
diam\,s_{i} \leq \frac{1}{10^{k_{0}}}\,diam\,\sigma_{i}\:;\;\forall i\in I,\forall j\in J\,;
\end{equation}

were $k_{0}$ is to be determined later.

\begin{figure}[h]
 \begin{center}
\includegraphics[scale=0.25]{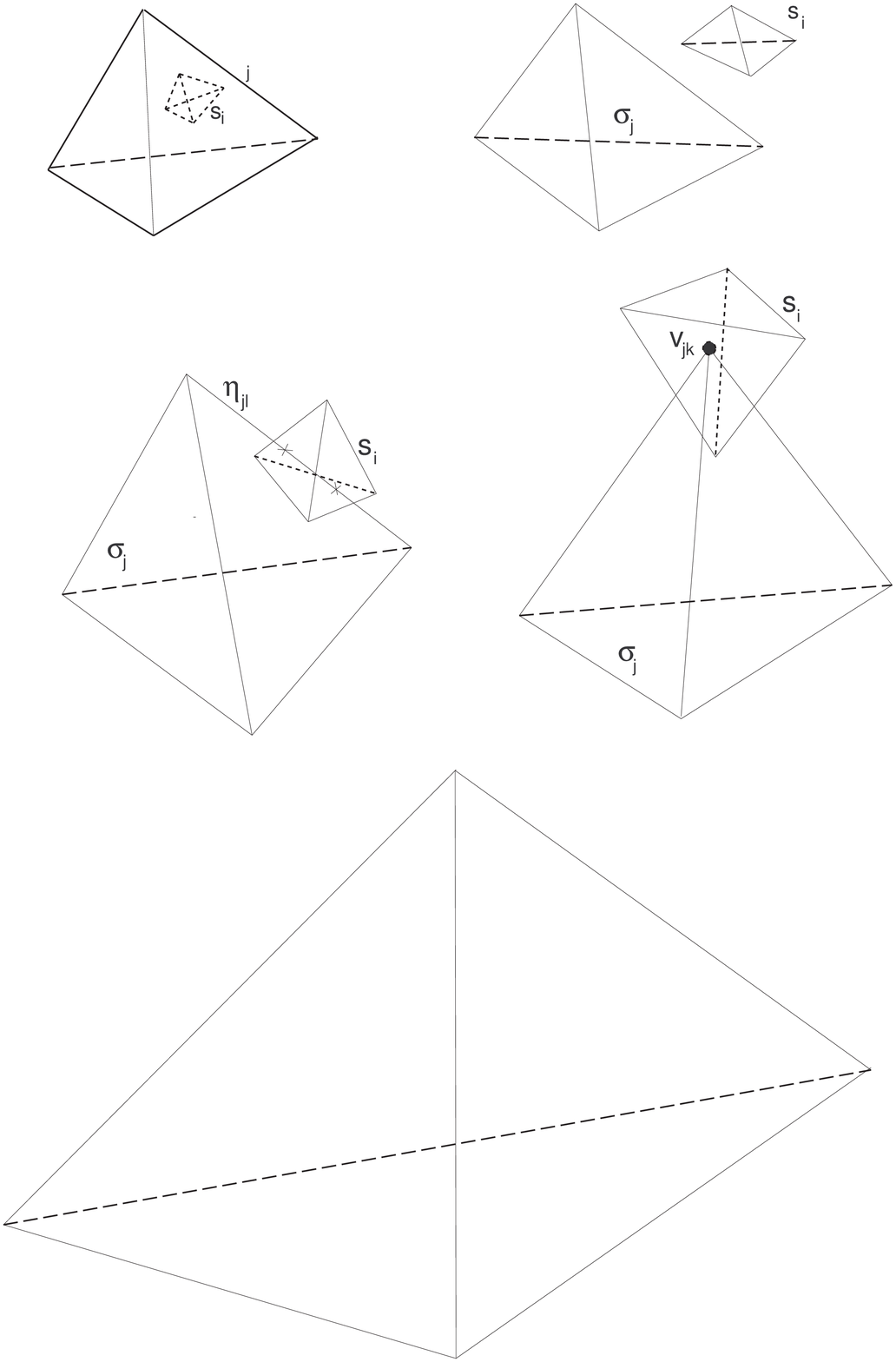}
\end{center}
 \caption{}
\end{figure}

By using general position arguments\footnote{\,completely analogous to the technique used in the proof of Theorem
1.3.} (see \cite{hu}\,, \cite{mun}\,),
 the relative positions of the $s_{i}$'s and the $\sigma_{j}$'s
 are now reduced to the following relevant possibilities:
\\ \hspace*{0.4cm} (a) $s_{i} \subset int \,\sigma_{j}$, or (a') $s_{i} \subset ext\,\sigma_{j}$\,;
\\ \hspace*{0.4cm} (b) $\exists \,\eta_{jl}\; s.t.\;  \eta_{jl}\, \cap \: int\,s_{i}\, \neq \,\emptyset$\,, where
                       $\eta_{jl}$ is an edge of $\sigma_{j}$, but $\exists\hspace{-0.2cm}/\; v_{im} \in \eta_{jl}$\,,
 \hspace*{0.47cm}      where $v_{im}$ is a vertex of $s_{i}$\,;
\\ \hspace*{0.4cm} (c) $\exists \, \nu_{jk}$ s.t. $\nu_{jk} \in int\,s_{i}$\,, where $\nu_{jk}$ is a vertex of $\sigma_{j}$\,;
\\ \hspace*{0.4cm} (d) $int\,s_{i} \,\cap\, int \,\sigma_{j} \neq \emptyset$\,, but we are not in one of the previous cases (see Fig. 3).

\subsection{Fattening $2$-dimensional Triangulations}
We start our triangulation "fattening" process by dealing with the $2$-dimensional case first:
\\Let $\sigma_{j_{0}} \in \Sigma$ be such that there exists a regular neighbourhood $N_{j_{0}}$ of $\sigma_{j}$,
triangulated by elements of $\mathcal{S} = \{ \sigma_{i}\}$.
 Now  $\mathcal{S}$ is partitioned by $\sigma_{j}$ into three disjoint families
 $\mathcal{S}_{0,1}\,,\: \mathcal{S}_{0,2}\,,\: \mathcal{S}_{0,3}$\,, where:
 \\ \hspace*{0.5cm}
  $\mathcal{S}_{0,1} = \{s_{i} \in \Sigma \;| \: s_{i} \subset int(\sigma_{j_{0}})\: {\rm or}\: s_{i} \subset int(\sigma_{j_{0}}) \}$
 \\ \hspace*{0.5cm}
  $\mathcal{S}_{0,2} = \{s_{i} \in \Sigma \;| \:\exists \:v_{0,k}\, - \,{\rm vertex\: of} \:\sigma_{j_{0}}\: {\rm s.t.}\: v_{0,k} \in int (s_{i}) \}$
 \\ \hspace*{0.5cm}
  $\mathcal{S}_{0,3} = \{s_{i} \in \Sigma \;| \: \exists!\:e_{0,l}\,-\,{\rm edge\:of}\:\sigma_{j_{0}}\: {\rm s.t.}\: e_{0,l}\, \cap \,int s_{i}\neq \emptyset \}$
\\It is easy to assure -- by eventual further subdivision and $\varepsilon$-moves\footnote{i.e. we "move a bit" the triangulation
$\mathcal{S} = \{s_{i}\}$ s.t. no intersection $s_{i} \cap \sigma_{j}$ is of one of the forbidden types. (See
\cite{mun} for definition and further details.)} -- that
\\ $\mathcal{S}_{0,3}  \cap  \mathcal{S}_{1,3} = \emptyset$\,; where $\mathcal{S}_{1,3}$ is the family corresponding to $\mathcal{S}_{0,3}$, induced by
$\sigma_{j_{1}}$, that is adjacent to $\sigma_{j_{0}}$.
\\ The intersections belonging to the family $\mathcal{S}_{0,2}$ are the principal generators of
"un-fatness", for $\angle(e_{0,l},e_{i,m})$ may be arbitrarily small, where $e_{i,m}\,; \\m = 1,2,3$ are the edges
of $s_{i}$ (see Fig. 3).

\begin{figure}[h]
 \begin{center}
\includegraphics[scale=0.3]{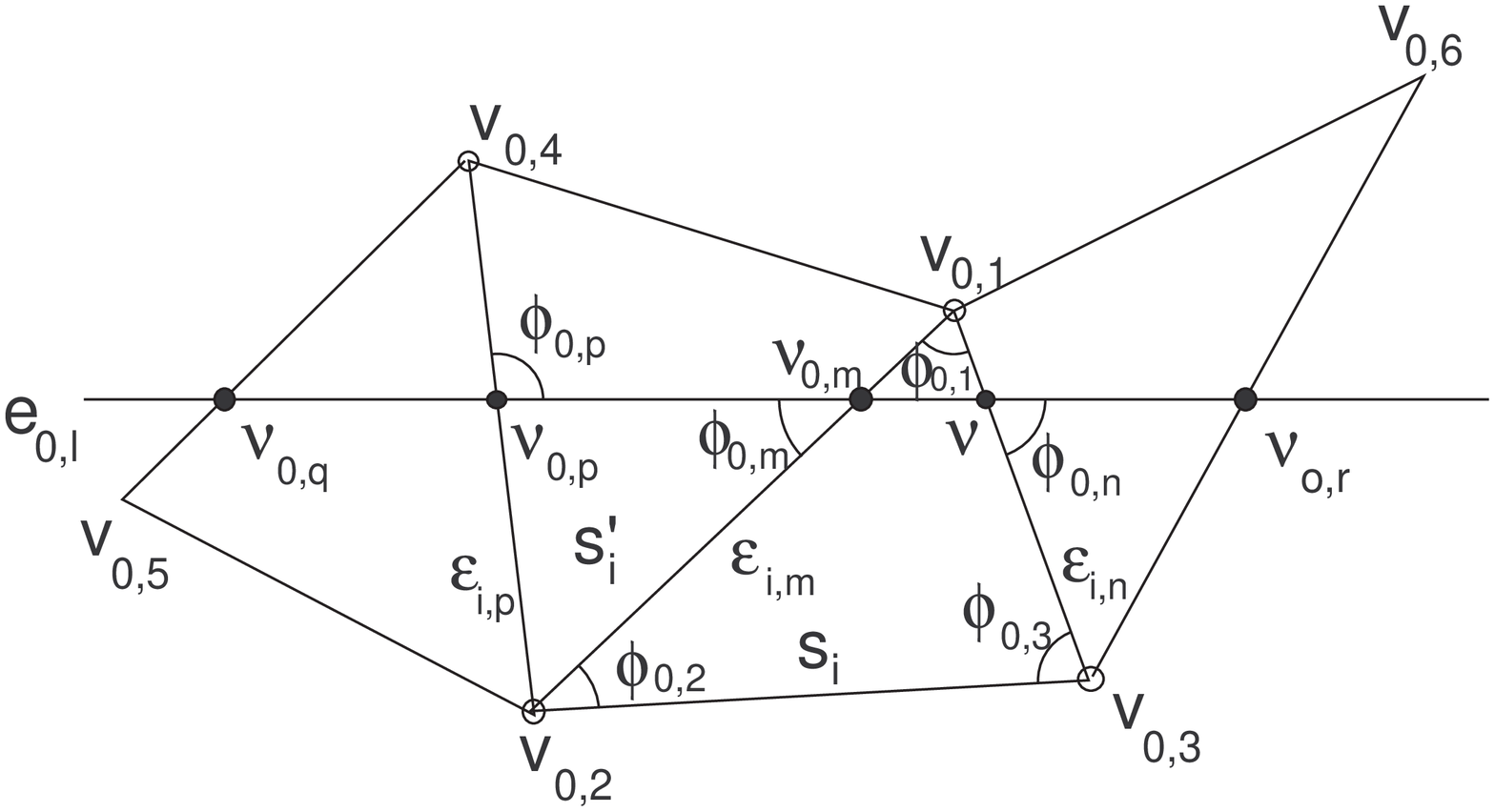}
\end{center}
 \caption{}
\end{figure}

Let $\nu_{0,p}\,,\;\nu_{0,m}\,,\;\nu_{0,n}$ be three consecutive intersection points of $e_{0,l}$ with edges of
two adjacent simplices $s_{i}$, $s'_{i}$ and let $\phi_{0,p}\,,\:\phi_{0,m}\,,\;\phi_{0,n}$ denote the resulting
angles (we always choose the acute angle) (see Fig. 3). Now it is not possible that two consecutive of the angles
$\phi_{0,p}\,,\:\phi_{0,m}\,,\;\phi_{0,n}$ are smaller than $\phi_{0}$: indeed, let us suppose that both
$\phi_{0,m} < \phi_{0}$ and $\phi_{0,n} < \phi_{0}$. Then $\phi_{0,1} > \pi - 2\phi_{0,2}$, so $\phi_{0,2} +
\phi_{0,3} < 2\phi_{0}$, and thus either $\phi_{0,2} < \phi_{0}$ or $\phi_{0,3} < \phi_{0}$, in contradiction to
the fatness of $s_{i}$\,.
\\ The fact above implies that we obtain two quadrilaterals which contain the "bad" points
$\nu_{0,p}$ and $\nu_{0,n}$ in their (respective) interiors, let them be: \nolinebreak[4] $Q_{1} = \Box
\nu_{0,1}\nu_{0,4}\nu_{0,5}\nu_{0,2}$ and $Q_{1} = \Box \nu_{0,2}\nu_{0,3}\nu_{0,6}\nu_{0,1}$ .
\\ We erase the segments  $\nu_{0,q}\nu_{0,p}\,,\: \nu_{0,2}\nu_{0,p}\,,\: \ldots \,,\:\nu_{0,n}\nu_{0,r}$
 and we replace them with segments that will "fattily" triangulate the quadrilaterals
 in question. These triangles have "big" angles (since their angles are belonging to fat triangles
 or are the sum of two such angles).
\\We distinguish between two cases:
\\(a) $Q_{i}$ is convex, and (b) $Q_{i}$ is not convex. We first take care of the simpler case.
\\ (a) Let $Q$ be a convex quadrangle such that all its angles are $\geq \phi_{0}$ and let $\overrightarrow{l}_{\!\!i,k}$ be the ray interior to $\angle A_{i-1}A_{i}A_{i+1}$
s.t. $\angle \overrightarrow{l}_{\!\!i,k}A_{i}A_{k} = \frac{1}{4}\alpha_{i}$\,; $i = 0,\ldots,3$\,; $k \in
\{i-1,i+1\}$ are considered $mod\,4$, of course. (See Fig. 4, where $Q^{*} = \chi_1\dots\chi_9$.) Then the rays
$\{\overrightarrow{l}_{\!\!i,k}\}_{i,k}$ generate a convex polygon $Q^{\ast} \subset int\,Q$. (See \cite{teza}.)

\begin{figure}[h]
 \begin{center}
\includegraphics[scale=0.3]{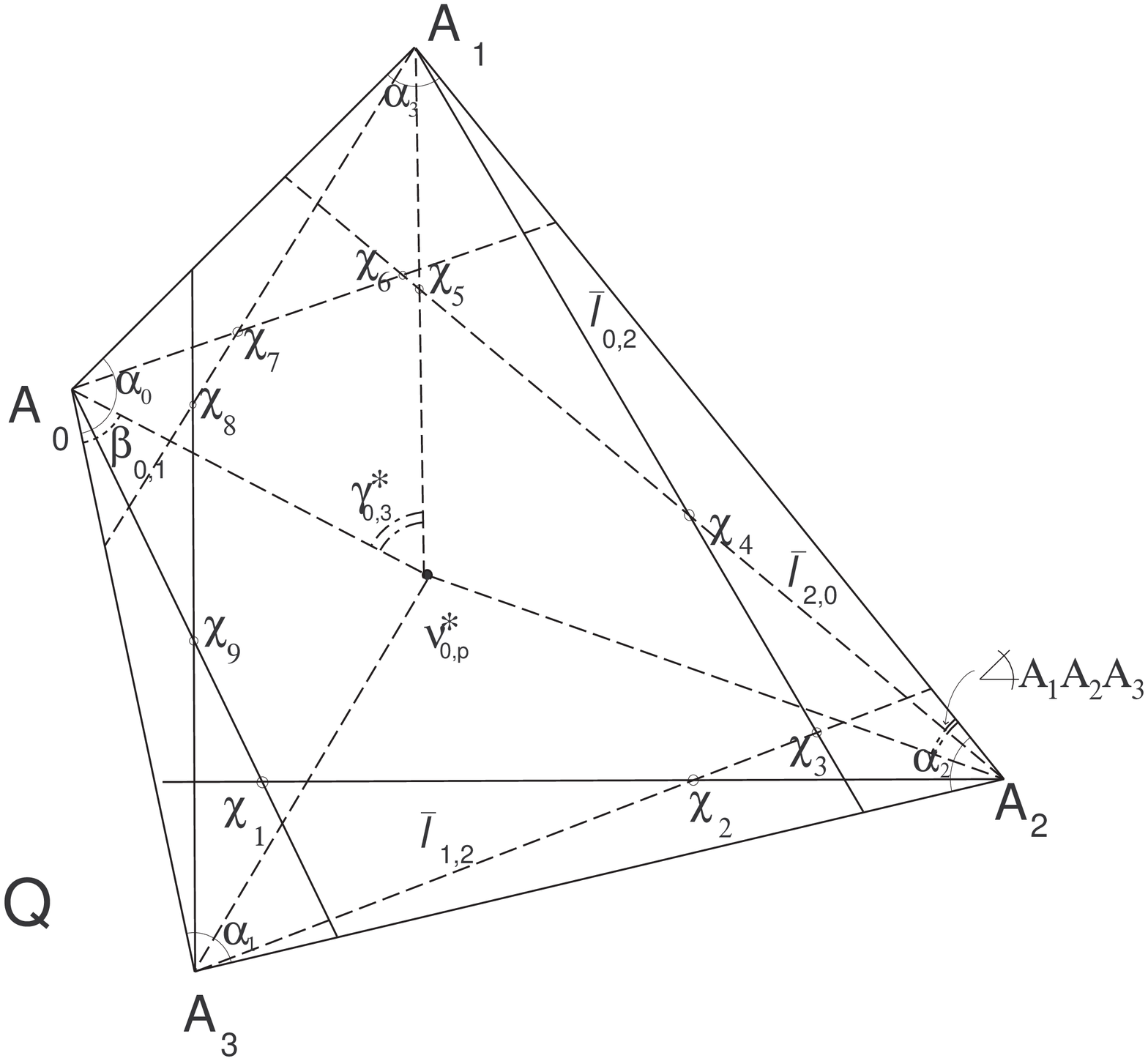}
\end{center}
 \caption{}
\end{figure}

By its very definition this quadrangle has the property that, for any $\nu^{\ast}_{0,p} \in int\,Q$ we have that:
\begin{equation}
\beta_{ik} = \angle \nu_{0,p}A_{i}A_{k} > \frac{\alpha_{i}}{4} \geq \frac{\phi_{0}}{4} \,;\; k \in
\{(i-1)mod\,4,(i+1)mod\,4\}.
\end{equation}
Also:
\begin{equation}
\gamma^{\ast}_{ik} = \angle A_{i}\nu^{\ast}_{0}A_{k} > \min \,\{\angle A_{i-1}A_{i}A_{i+1},\angle
A_{i}A_{i+1}A_{i-1}\}\; .
\end{equation}
(See Fig. 5) But one one of the angles $\angle A_{i-1}A_{i}A_{i+1}$ and $\angle A_{i}A_{i+2}A_{i-1}\}$ belongs to
one of the original fat triangles, so:
\begin{equation}
\phi^{\ast}_{ik} > \phi_{0} \;.
\end{equation}
So, from $(4.1.)$ and $(4.3.)$ it follows that the triangles $\triangle A_{i}\nu^{\ast}_{0}A_{i+1}\,;\; i =
0,\dots,3\,(mod4)$ are fat.

\begin{figure}[h]
\begin{center}
\includegraphics[scale=0.25]{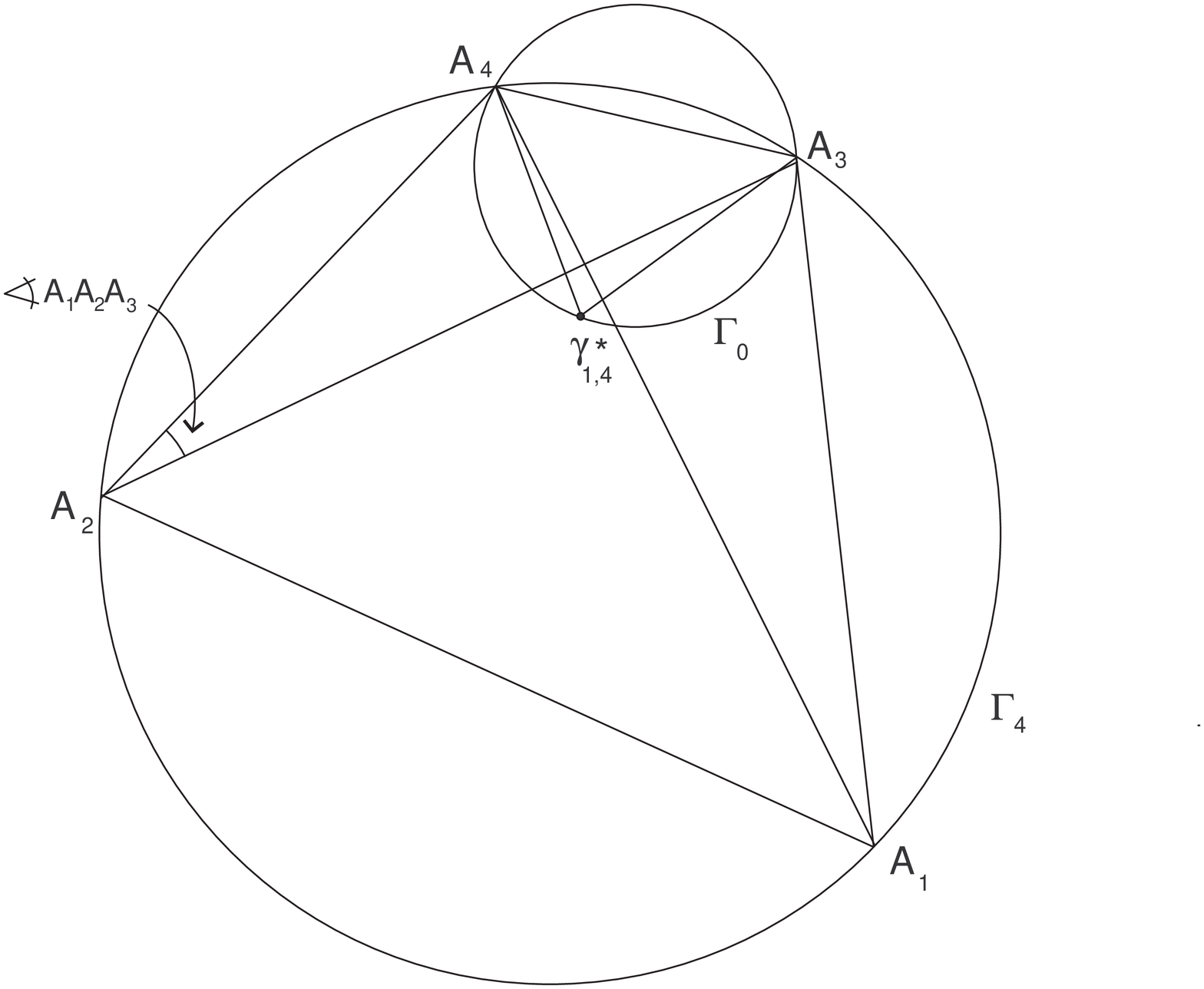}
\end{center}
\caption{}
\end{figure}

(b) In this case, rather than tracking back our steps through the same argument as in the previous case; we prefer
to dissect $Q$ into two triangles and one convex quadrilateral, in the following way: if $\Box
A_{0}A_{0}A_{0}A_{3}$ is such that $\angle A_{3}A_{0}A_{1} > \pi$ and such that $\nu_{0,p} \in A_{2}$\,, then
consider the bisectors $A_{0}B_{12}$ of and $\angle A_{1}A_{0}A_{2}$ and $A_{0}B_{23}$ of $\angle
A_{2}A_{0}A_{3}$. (See Fig. 6.)

\begin{figure}[h]
\begin{center}
\includegraphics[scale=0.25]{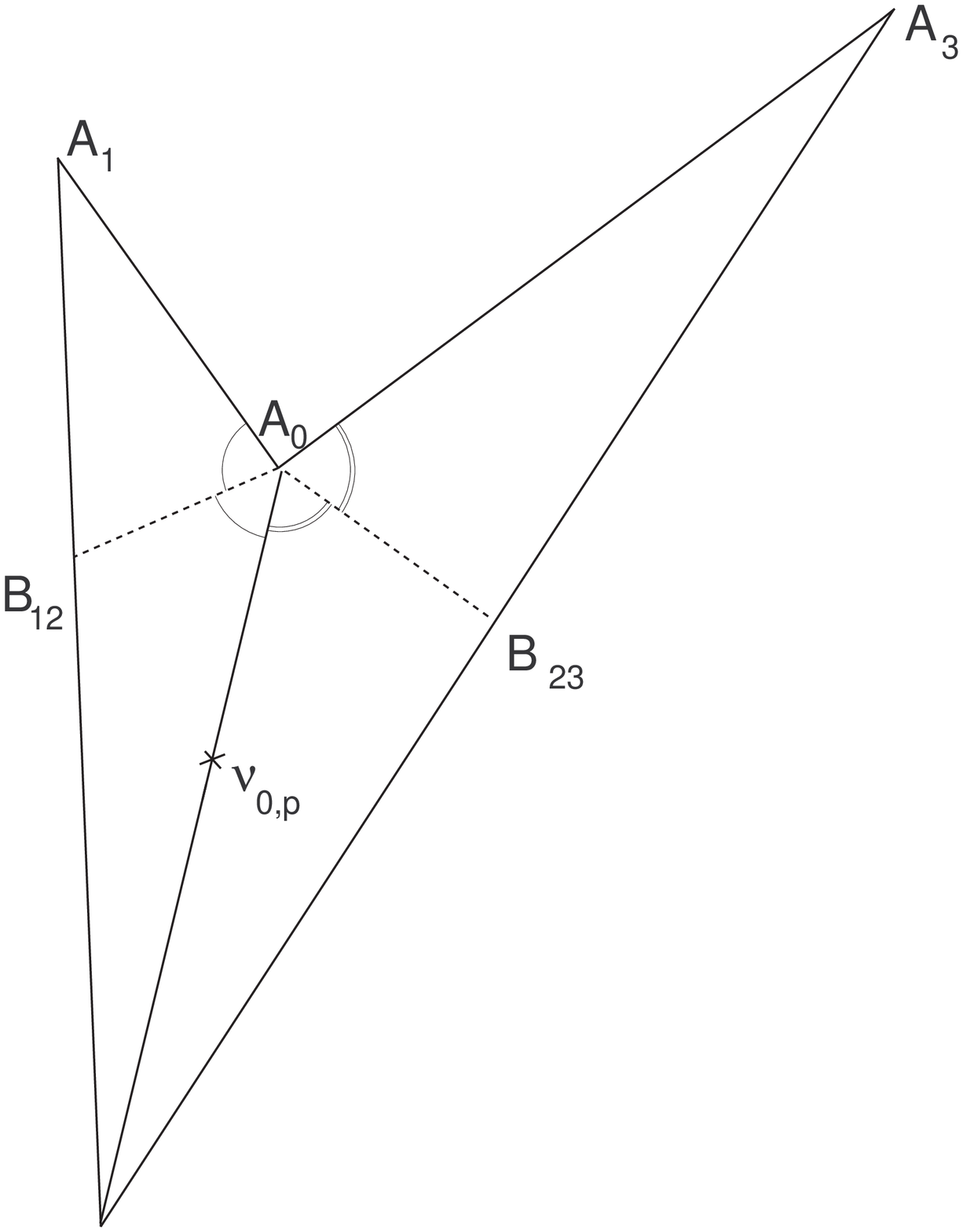}
\end{center}
\caption{}
\end{figure}

Now, instead of returning once over to the argument used in case (a), it is easy to proceed directly and show that
the triangles  $\triangle A_{1}A_{0}B_{12}\,,\ldots,\, \triangle A_{0}B_{23}A_{3}$  are fat. (See \cite{teza} for
details.)
\\We are now faced with a new fat triangulation.
However, new points of intersection with $e_{0,l}$ are introduced: let them be $\nu_{0,p-1}$ and $\nu_{0,p+1}$. If
\[e_{0,l} \cap \Box A_{0}B_{12}A_{2}B_{23} \cap \{B_{12}B_{23}\} = \emptyset\,,\]
then the same argument we used for the original triangulation shows that $\nu_{0,p-1}$ and $\nu_{0,p+1}$ are $\geq
\phi_{0}$. If
\[e_{0,l} \cap \Box A_{0}B_{12}A_{2}B_{23} \cap \{B_{12}B_{23}\} \neq \emptyset\,,\]
we can employ either one of the following two methods to remedy the situation: (a) use the general position
technique again and bring the new triangulation to the required position; or (b) consider, instead of the bisector
$A_{0}B_{23}$ the meridian $A_{0}M_{23}$ ($M_{23} \in A_{2}A_{3}$) -- see \cite{teza} for details.

We still have to face the problem of "mashing" the triangulations "over the wave front" of the "$s_{i}$"-s. Away
from the vertices of the complex $\mathcal{S}$, we are faced with two possibilities: $\mathcal{S}$ contains: (a)
one or (b) two of the vertices of $\sigma \in \Sigma$. We shall deal first with:
\\ {\sf Case (b)} By further dividing\footnote{\,since we are in the planar case, the division may be done by lines parallel to the edges of $s$.}
the triangles of:
\begin{equation}
Front\,\mathcal{S} = \{s \in \mathcal{S}\, |\, \exists \,e - edge\: of\: s\: {\rm s.t.}\: e \in \partial
\mathcal{S} \} \;,
\end{equation}
we are able to "erase" the families $\mathcal{S}_{1}$ and $\mathcal{S}_{2}$\,, where $\mathcal{S}_{1}\: \cup \:
\mathcal{S}_{2} = Front\,\mathcal{S}\: \cap \: \{e_{i}\}$, $i = 1,2$ (see Fig. 7 (a),\:(c)\,) our only problem
being that some of the simplices $s \in \mathcal{S}$ may intersect the edges $e_{i}$ of $\sigma$ at an angle
$\phi_{s,l_{i}} < \phi_{0}\,,\; i=1,2$. However, the angles $\phi^{2}_{s,l_{i}}$ (see Fig. 7 (b)) are,
 by the previous argument  $\geq \phi_{0}$.

\begin{figure}[h]
\begin{center}
\includegraphics[scale=0.25]{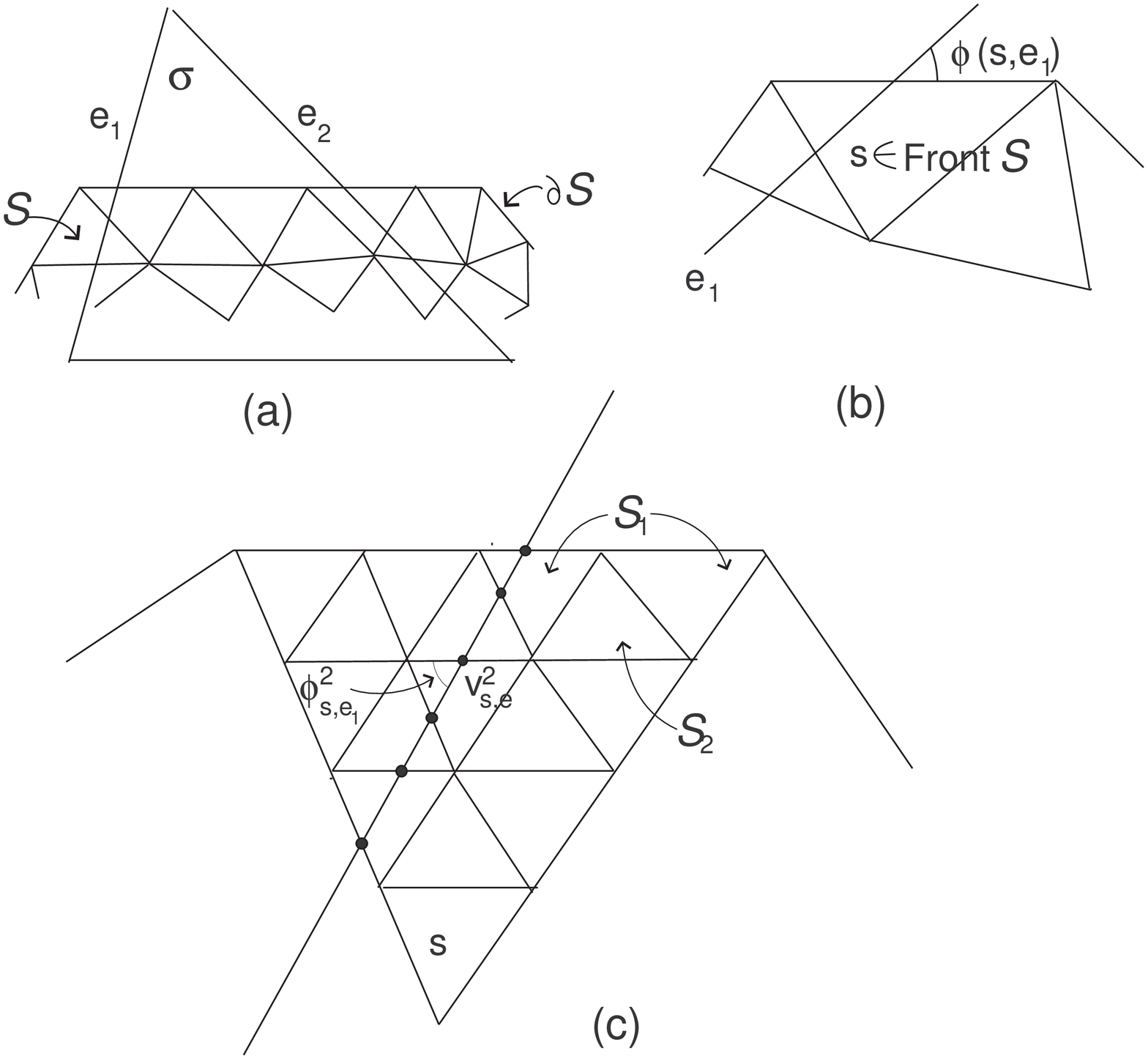}
\end{center}
\caption{}
\end{figure}

Let $\mathcal{S}' = \mathcal{S} \setminus Front\,\mathcal{S}$. Then $\partial \mathcal{S}'$ will still be convex,
so we can consider the {\it joins} ({\it cones}) $J(v_{12},\varepsilon_{k})$, where $\varepsilon_{k}$ are the
edges of $\partial \mathcal{S}'$ that are included in $\sigma$, and $v_{12} = e_{1} \cap e_{2}$ (see Fig. 8).

\begin{figure}[h]
\begin{center}
\includegraphics[scale=0.25]{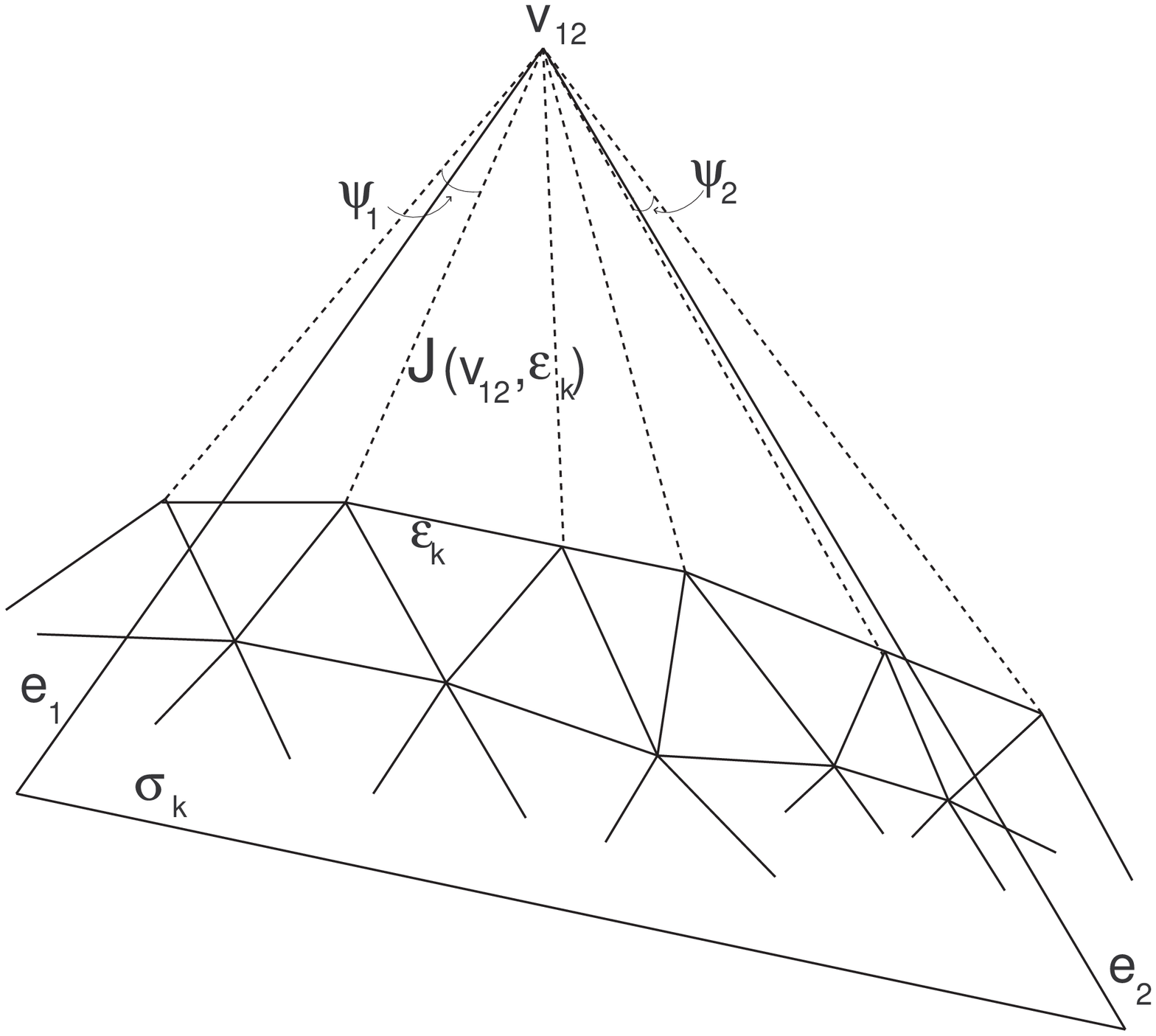}
\end{center}
\caption{}
\end{figure}

Now, since the number $n_{0}$ of conjugacy classes of elliptics is finite, we can consider
\begin{equation}
\delta_{0} = min\{\delta_{1},\delta_{2},\ldots,\delta_{n_{0}}\}\:;
\end{equation}
where the bounds $\delta_{i}\,,\; 1 \leq i \leq n_{0}$ are those given by Theorem 2.1. Then in $\partial M_c
\times [3/4,1)$ (or alternatively in the tubes $T_{0,i} = N_{i}(\delta_{0}/4) \setminus N_{i}(3\delta_{0}/16)$)\,,
the radii $r_{0}$, $R_{0}$ of the simplices $\sigma_{0}$, \nolinebreak[4]$\sigma_{0} \, \bigcup
\raisebox{-0.7em}{\hspace{-0.5cm}\tiny $i\in N$}J_{0} \neq \emptyset$, where $J_{0} = \bigcup T_{0,i}$\,, will be
uniformly bounded.\footnote{\,See \cite{pe}.}
\\Thus, there exist numbers $m_{0}$ and $m_{1}$ such that:
\begin{equation}
m_{0} \leq diam(\sigma) \leq m_{1} \,;\: \forall \,\sigma \: {\rm s.t.}\: \bar{\sigma} \cap T_{0} \neq \emptyset
\,.
\end{equation}
Therefore, exists $k_{1} \in N_{+}$ such that:
\begin{equation}
\frac{1}{10^{k_{1}}}diam(\sigma) \, \leq \, diam(s) \,; \: \forall \,s,\sigma \: s.t. \: \bar{s},\,\bar{\sigma}
\cap T_{0} \neq \emptyset \,.
\end{equation}
This, in conjunction with (4.2), assures us that the number $\rho_{1}$  of triangles $s_{i}$ that intersect the
edge $e_{0,l}$ (see Fig. 9.\,(a)) will be bounded by two natural numbers, i.e. $\exists \: n^{1}, n^{2} \in N$
s.t.
\begin{equation}
n^{1} \leq \rho_{1} \leq  n^{2}\, .
\end{equation}

\begin{figure}[h]
\begin{center}
\includegraphics[scale=0.3]{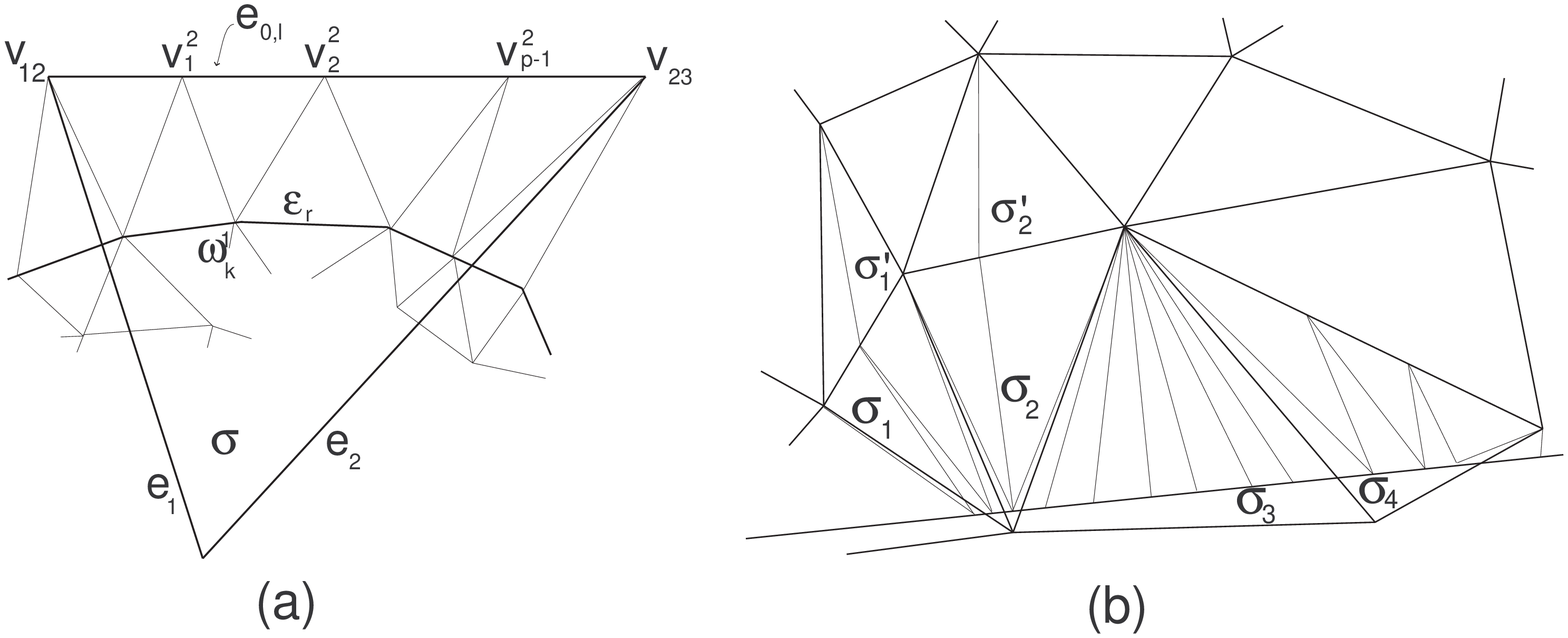}
\end{center}
\caption{}
\end{figure}

Moreover, our last subdivision ensures us that the number $\rho_{2}$ of triangles  intersecting $\partial
\mathcal{S}'$ is:
\begin{equation}
\rho_{2} \leq 3(\rho_{1}-1)+2\;,
\end{equation}
so there exists $\lambda_{0} > 0$ s.t. the angles $\psi_{k} = \angle v_{12}, \psi^{1}, \psi^{2}$ (see Fig. 8)
satisfy
\begin{equation}
 \psi_{k} \geq \lambda_{0}\varphi_{0}\,, \; \psi^{i} \geq \lambda_{0}\varphi_{0}\,, \; i \,=\, 1,2\,;
\end{equation}
as desired.
\\ We still have to deal with
\\ {\sf Case (a)}  We repeat once more the procedure used for the subdivision of $\partial \mathcal{S}$
and $Front \,\mathcal{S}$ that we employed in Case (b). We divide the edge $e_{2}$ into $\rho_{2} - 1$ equal
segments and consider the joins (cones) $J(v^{2}_{j},\varepsilon_{j+1})$ and $J(w_{k}',v^{2}_{k}v^{2}_{k+1})$ (see
Fig. 9\,(a)). Using the same arguments as before one easily checks that the resulting triangulation of $s
\,\backslash \,\mathcal{S}$ will be fat. Moreover, although this procedure dramatically reduces the "fatness" of
the next stratum of simplices of the family $\Sigma$, it leaves the other strata unchanged.
\\This procedure takes care of the intersections of $\Sigma$ with "the front wave" \\$\mathcal{S} \cap \{x\,|\,d_{hyp}(A_i,x) = \delta_{0}/4\}$. To deal with the case $\mathcal{S} \cap \{x\,|\,d_{hyp}(A_i,x)  =
3\delta_{0}/16\}$ one proceeds along the same lines and then fits the new triangulation $\mathcal{S}'$ to
$\mathcal{S}$ in a properly chosen tubular region $J_{1}$, using, instead of $\sigma$, triangles $s_{0} \in
\mathcal{S} \cap \mathcal{T}_{1}$.
\\ We do have yet to contend with the problem posed by "corners" i.e. by triangles $s \in \mathcal{S}$
such that (i) $s \in Front\,\mathcal{S}$ and (ii) $\partial s  \cap  \partial \mathcal{S} = {v}$ (were $v$ is a
vertex). The two cases that ensue may be treated with the methods developed before. (See \cite{teza}\,.\.)
\begin{rem}
We can dispense with theses last considerations altogether by considering the following: we are concerned -- in
fact -- only with patching together triangulations already contained in a $\delta_{0}$-neighbourhood of an axis.
But the geometric triangulations employed partitioned these neighbourhoods into "levels" or "heights", so one can
fit the triangulation of two consecutive levels -- say "$m$" and "$m+1$" -- by considering the intermediary
adjusting patch delimited by the levels  "$m - \frac{1}{2}$" and  "$m + \frac{1}{2}$". (A finite number of further
barycentric subdivisions may still be required.)
\end{rem}
This concludes the "fattening" process for the two dimensional triangulations.

\subsection{Fattening $3$-dimensional Triangulations}
We shall divide the "fattening" process of two intersecting $3$-dimensional tetrahedra into two parts: A) The
"fattening" of the $2$-dimensional intersection between the triangles belonging to $\partial \mathcal{S}$ and a
face $f_{123} = \triangle v_{1}v_{2}v_{3}$ of a tetrahedron $\sigma \in \Sigma$\,, and B) The extension of the new
triangulation to a fat $3$-dimensional triangulation.

\subsubsection{Fattening $2$-dimensional intersections}
 Let us start by noticing that, since the fatness of the tetrahedra $s_{i} \in \mathcal{S}$ is bounded from
below, so will be their dihedral angles\footnote{\, by the "trihedral sinus formula" (see \cite{bac})}.  It
follows that even after the partition of $f_{123}$ into triangles $\{\tau_{k}\}$ and quadrilaterals $\{
\eta_{j}\}$, to a triangulation $\{ \tilde{\tau}_{k}\}$, the number of triangles around each vertex will be
bounded by a natural number $m_{1}$. We shall exploit this fact to our advantage, so that we will be able to
replace $\{\tilde{\tau}_{k}\}$ by a fat triangulation $\{ \tilde{\tau}_{k}^{\ast}\}$ that is "close" to
$\{\tilde{\tau}_{k}\}$.
 \\Indeed, if we denote by $\tilde{\alpha}^{0}_{j}\,,\:j = 1,\dots,\bar{m}_{0}, \:\bar{m}_{0} \leq m$\,;
 the angles of the triangles $\{ \tilde{\tau}_{k}^{0}\}$ that are respectively adjacent
to the vertex $\tilde{\nu}^{0} = \bigcap_{j} \tilde{\tau}_{j}^{'}$  (see Fig. 10), then there are two sources of
"thiness": either $\tilde{\alpha}^{0}_{j}\:=\:\alpha^{0}_{j}$ is smaller than $\varphi_{0}$, for some $j_{0} \in
\{1,\ldots,\bar{m}\}$\; or one of the angles $\tilde{\alpha}^{0}_{j_{1}}\,\tilde{\alpha}^{0}_{j_{1}+1}$ produced
by the division of the quadrilateral $\eta_{_{1}}$ into two triangles: $\tilde{\tau}^{0}_{j_{1}}$ and
$\tilde{\tau}^{0}_{j_{1}+1}$.

\begin{figure}[h]
\begin{center}
\includegraphics[scale=0.33]{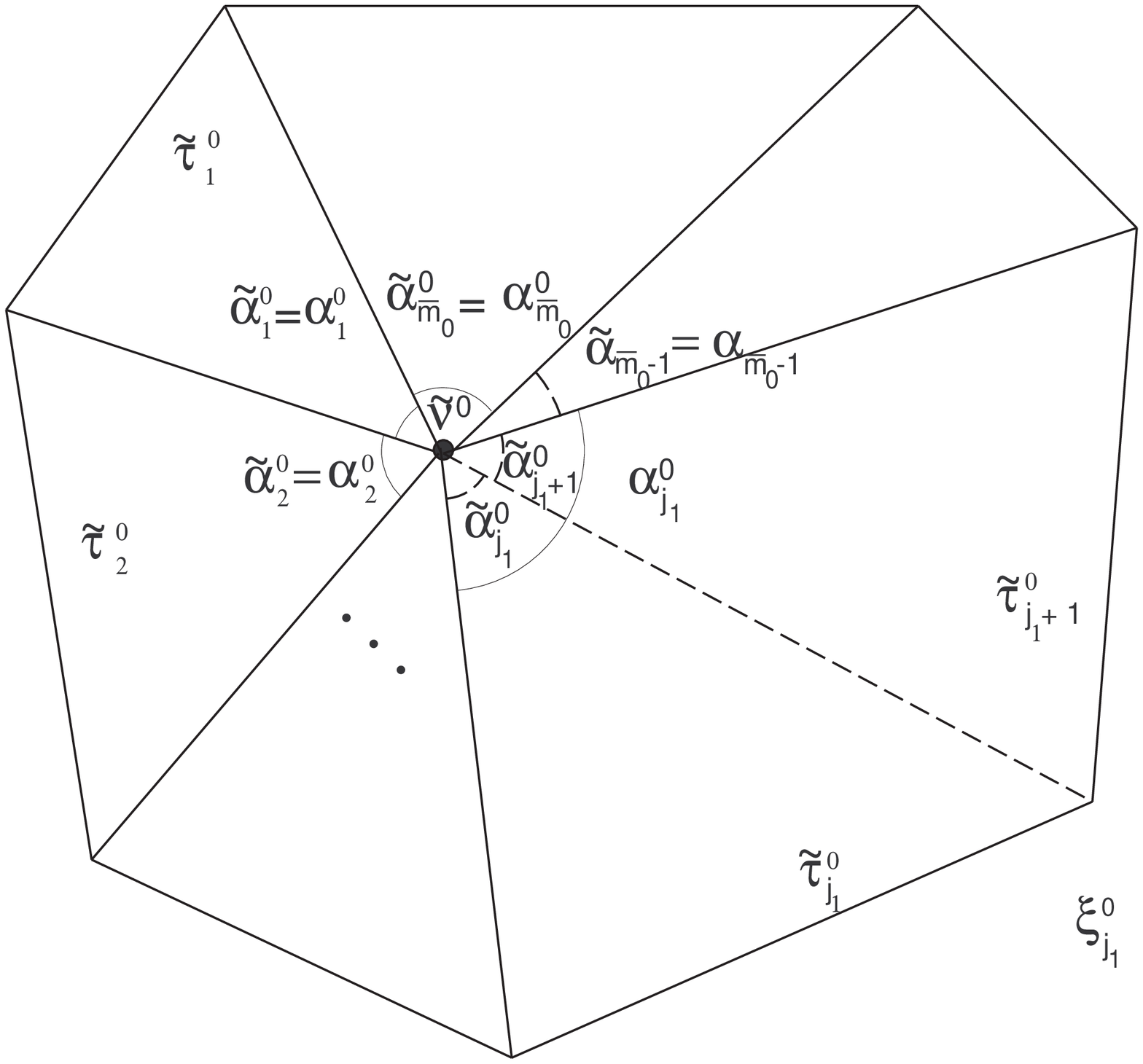}
\end{center}
\caption{}
\end{figure}

If we bisect the angles $\tilde{\alpha}^{0}_{j}\,,\;j=1,\ldots,\bar{m_{0}}$\,; we will receive angles
\nolinebreak[4] $\tilde{\alpha}^{0}_{jk}\,,\\ \:j\:=\:1,\ldots,\bar{m_{0}}$\,; $k = 1,2$. Let us consider the
angles $\tilde{\beta}^{0}_{jk}\,,\:j = 1,\ldots,\bar{m_{0}}$\,; $k = 1,2$\,; where: $\tilde{\beta}^{0}_{j1} =
\frac{\tilde{\alpha}^{0}_{j-1} + \tilde{\alpha}^{0}_{j}}{2}$\,,\; $\tilde{\beta}^{0}_{j2} =
\frac{\tilde{\alpha}^{0}_{j} + \tilde{\alpha}^{0}_{j+1}}{2}$\,.\footnote{\,The indices are to be taken $mod(\bar
m_{0})$, of course.} If each end every of the angles $\tilde{\beta}^{0}_{jk}$ defined above is greater than
$\varphi_{1} = \varphi_{0}/10$, then we desist and proceed towards part. But it may be that, for instance, both
$\tilde{\alpha}^{0}_{j-1}$ and  $\tilde{\alpha}^{0}_{j}$ are smaller than $\varphi_{0}/5$. If this happens we
continue the process of "mixing the angles".  To be more precise: let us delate -- for commodity reasons -- the
upper index  "$0$" in the enumeration of the angles "$\beta$", and denote them by $\tilde{\beta}_{j}, \; j\,
=\,1,\dots,\bar{m}_{0}$\,; and let us form the sequence of angles $\tilde{\beta}_{j}'$, $\tilde{\beta}_{j}''$,
$\tilde{\beta}_{j}'''$, etc.\,, where $\tilde{\beta}_{j}' = \frac{1}{2}( \tilde{\beta}_{j-1} +
\tilde{\beta}_{j+1})$, $\tilde{\beta}_{j}'' = \frac{1}{2}( \tilde{\beta}_{j-1}' + \tilde{\beta}_{j}')$, and so on.
But this process will halt -- inasmuch as we are concerned -- in a finite number of steps, for the following
inequality holds:

\begin{equation}
\tilde{\beta}^{\bar{m}_{0}}_{j} > \frac{\alpha_{1}+ \cdots +\alpha_{\bar{m}_{0}}}{2^{\bar{m}_{0}+1}} \,; \;
j\,=\,1,\ldots,\bar{m}_{0}\;;
\end{equation}
that is:
\begin{equation}
\beta^{\bar{m}_{0}+1}_{j} > \frac{2\pi}{2^{\bar{m}_{0}+1}} > \frac{\varphi_{0}}{10^{k}} \,; \;
j\,=\,1,\ldots,\bar{m}_{0}\;;
\end{equation}
and so:
\begin{equation}
\beta^{\bar{m}_{0}-1}_{j} > \varphi_{1} =  \frac{\varphi_{0}}{10^{\bar{k}}} \,;\;1,\ldots,\bar{m}_{0}\;;
\end{equation}
where $\bar{k}$ is the least natural power that satisfies the right-handed inequality in (4.12).
\\We shall use the bound above in order to produce a fat triangulation. However, some care is needed in doing this, for in general, both the
number of iterations used for each vertex and the number of bisectors $b^{k}_{j}\,,\; j =
1,\ldots,\bar{m}_{0}\,;\; k = 1,2,3$\,; that intersect the triangle $\tilde{\tau}$ will be different, thus
affecting the sizes of the angles $\tilde{\alpha}^{k}_{j}$.

\begin{figure}[h]
\begin{center}
\includegraphics[scale=0.33]{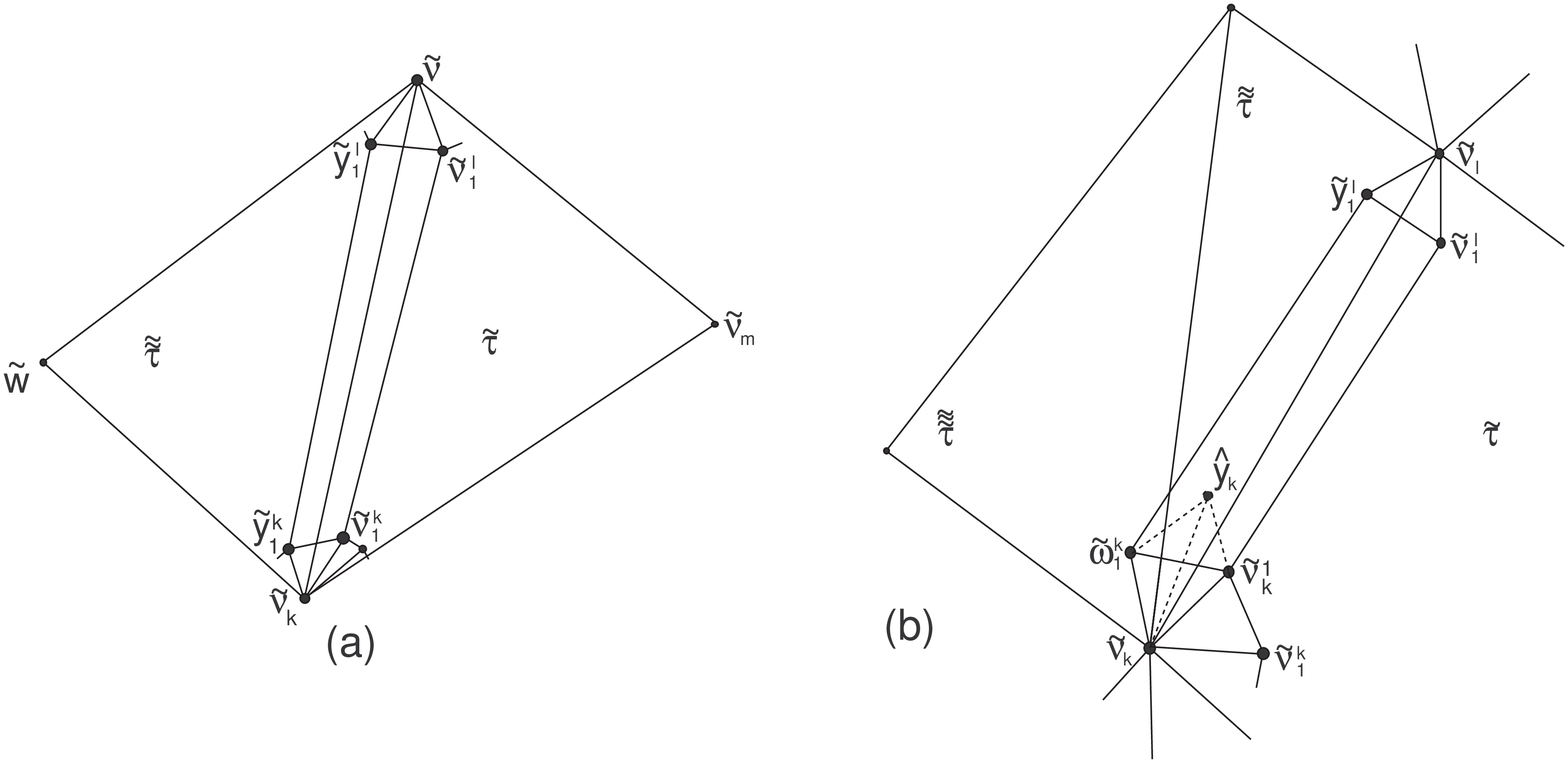}
\end{center}
\caption{}
\end{figure}

Let $b$ denote the shortest of the segments $b^{k}_{j} \cap \tilde{\tau}\,,\; k = 1,2,3.$  (By elementary
geometry, $b$ should be the segment $b^{k_{0}}_{i_{0}}$ that is the nearest to the shortest side of
$\tilde{\tau}$, where $\tilde{\alpha}^{k_{0}}$ is the smallest angle of $\tilde{\tau}$.)\footnote{The best
possibility occurs, of course, when we have to use the bisection only once, so we will have only (ordinary)
bisectors $b^{1},\,b^{2},\,b^{3}$ that will meet, of course, at the barycenter $\tilde{b}$ of $\tilde{\tau}$.}
 Moreover, let us consider segments
 $\tilde{\nu}_{k}\nu^{k}_{i} \subset b^{k}_{i}$\,, such that $length(\tilde{\nu}_{k}\nu^{k}_{i}) = b/10^{k_{2}} = \hat{b}$
where $k_{2}$ is chosen in such a manner that the quadrilateral $\Box \nu^{k}_{1}\nu^{l}_{1}y^{k}_{1}y^{l}_{1}$ is
simple; where $\nu^{k}_{1},\, \nu^{l}_{1}$ denote the vertices  closest to the edge $\nu_{k}\nu_{l}$ and
$y^{k}_{1},\, y^{l}_{1}$ play the same role in the adjacent triangle $\tilde{\tilde{\tau}}$ (see Fig. 11 (a)). It
may be that we will have to consider  $\Box \nu^{k}_{1}\nu^{l}_{1}\tilde{y}^{k}_{1}\tilde{y}^{l}_{1}$ (or any
permutation of indices) -- see Fig. 11 (b). (Here $\tilde{w}^{k}_{1}$ plays  in  the  triangle
$\tilde{\tilde{\tilde{\tau}}}$ the role $\tilde{\nu}^{k}_{1}$ plays in the triangle $\tilde{\tau}$.) In this case
consider a point $\hat{y}_{k} \in int\,\tilde{\tilde{\tilde{\tau}}}$, such that:
\\$(i) \;\;\;\;length(\hat{y}_{k}) = \hat{b}$\,;
\\$(ii_{1})\, \angle \tilde{w}^{k}_{1}\hat{\nu}_{k}\hat{y}_{k}\; \geq \; \frac{\varphi_{1}}{2}\;=\;\varphi_{2}$\,;
\\$(ii_{2})\, \angle \hat{y}_{k}\hat{\nu}_{k}\hat{\nu}^{k}_{1}\; \geq \; \frac{\varphi_{1}}{2}\;=\;\varphi_{2}$\,.

\begin{rem} The existence of the positive integer $k_{2}$ with the desired proprieties is guaranteed by the
fact that the angles $\tilde{\alpha}^{k},\, \tilde{\alpha}^{l},\,\tilde{\alpha}^{k}_{i},\,\tilde{\alpha}^{l}_{i},
\tilde{\alpha}^{m}$ and   $\tilde{\tilde{\alpha}}^{m}$ are  $\geq \varphi_{1}$, thus
\[ c_{1}b \leq length(\tilde{\nu}_{k}\tilde{\nu}_{k}) \leq  c_{2}b\,; \]
\[\;\, c_{1}'b' \leq length(\tilde{\nu}_{k}\tilde{\nu}_{k}) \leq  c_{2}'b'\,. \]
\end{rem}

In consequence we are facing the situation depicted in Fig. 12\,, where we also illustrated the conning  of the
segments $\tilde{\nu}^{k}_{i}\tilde{\nu}^{k}_{i+1}$ and $\tilde{\nu}^{k}_{1}\tilde{\nu}^{l}_{1}$ from the
barycenter $\tilde{\tilde{b}}$ of $\tau$; $k = 1,2,3;\; l = 1,2,3$.

\begin{figure}[h]
\begin{center}
\includegraphics[scale=0.38]{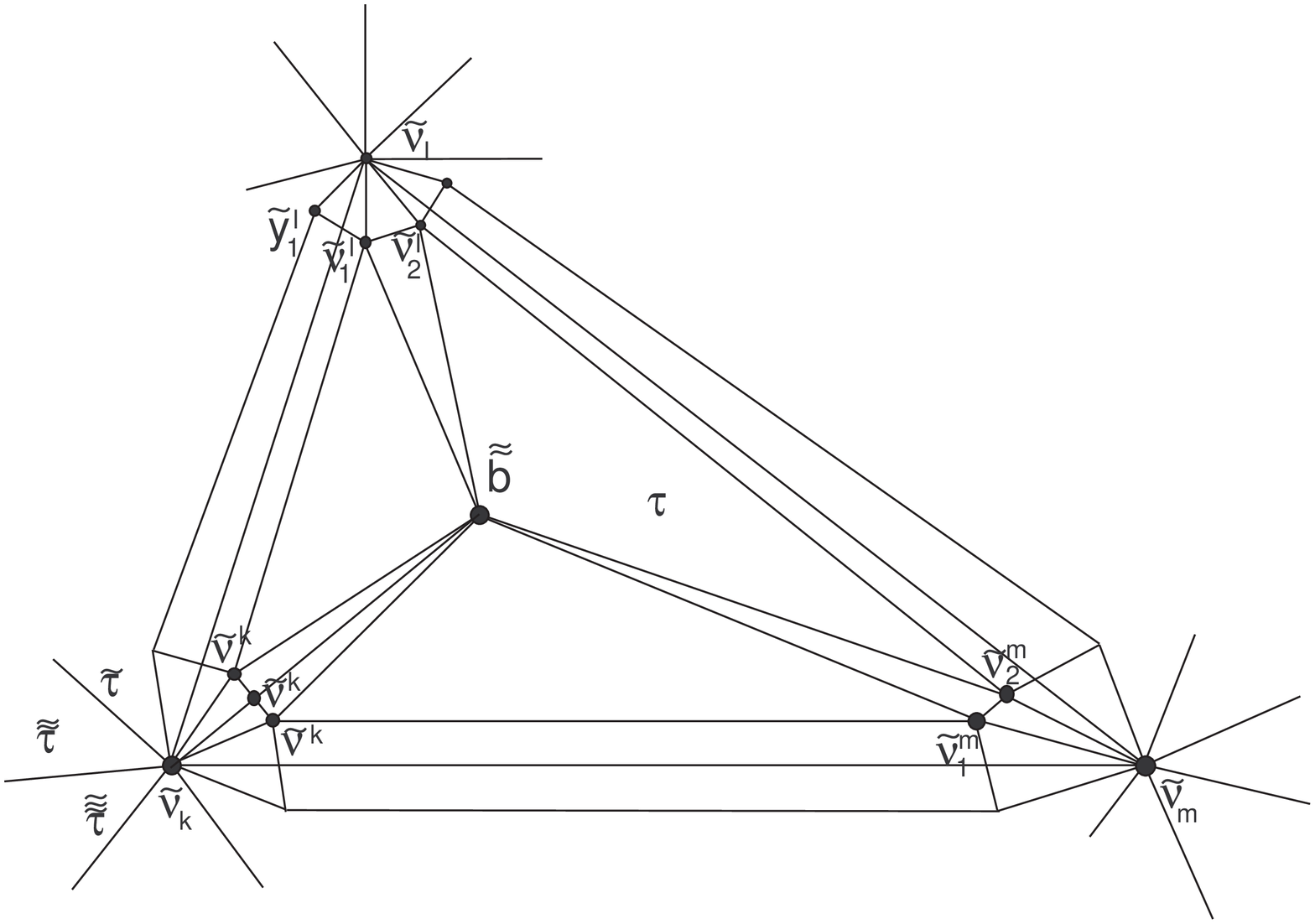}
\end{center}
\caption{}
\end{figure}

Now we have to show the fatness of three types of polygons:
\\ \hspace*{0.5cm} 1) triangles $T_{ij}\, = \,\triangle \tilde{\nu}^{k}_{i}\tilde{\nu}^{l}_{j}\tilde{\tilde{b}}$\,;
\\ \hspace*{0.5cm} 2) quadrilaterals $Q_{kl}\, =\, \Box \tilde{\nu}^{k}_{1}\tilde{\nu}^{l}_{1}\tilde{y}^{l}_{1}\tilde{y}^{k}_{1}$\,;
\\ \hspace*{0.5cm} 3) triangles $T^{kl}\, =\, \triangle \tilde{\nu}^{k}_{1}\tilde{\nu}^{l}_{1}\tilde{\tilde{b}}\,,\;  k,l\,=\,1,2,3$\,.
\\ We start by noticing that "fatness" of a quadrilateral means that:
\\ \hspace*{0.5cm} a) its angles are bounded from below, and:
\\ \hspace*{0.5cm} b) the ratios $l_{\lambda}l_{\iota}\,,\;\lambda,\iota = 1,\ldots,4m$ between the lengths of its sides are also bounded from below.

 An easy computation (see \cite{teza}\,) shows that
\begin{equation}
   length(\tilde{y}^{l}_{1}\tilde{\nu}^{l}_{1}) \; \geq \; 2\hat{b}\sin{\frac{\pi}{\bar{m}_{0}}}\;;
\end{equation}

\begin{equation}
   length(\tilde{y}^{l}_{1}\tilde{\nu}^{k}_{1}) \; \geq \; length(\tilde{\nu}^{k}_{1}\tilde{\nu}^{l}_{1}) - 2\hat{b}\,;
\end{equation}

\begin{equation}
   \angle \tilde{\nu}^{l}_{1}\;\geq\; \arctan{\frac{\hat{b}\cos{\frac{2\pi}{\hat{m}_{0}}}}{length(\tilde{\nu}^{k}\tilde{\nu}^{l}) - 2\hat{b}\sin{\frac{2\pi}{\bar{m}_{0}}}}}\;;
\end{equation}

where $\hat{b}\;=\;length(\tilde{\nu}^{k}\tilde{y}^{k}_{1})$\,.\footnote{\, and similar formulas hold for the
other pairs of vertices.}
\\ Therefore, the arguments employed before show that $\Box Q_{kl}$ is decomposable into fat
triangles, so case 2) is dealt with.
\\ In a manner similar to that used in case 2) show that the triangles of types $T_{ij}$ and
$T_{kl}$ are also fat; indeed:

\begin{equation}
   \angle \tilde{\nu}^{l}_{i}\tilde{\tilde{b}}\tilde{\nu}^{k}_{i+1} \; < \;2\arctan{\frac{3\pi}{2\bar{m}_{0}b^{l}}}\;;
\end{equation}
where

\begin{equation}
     lenght(\tilde{y}^{l}_{1}\tilde{\nu}^{k}_{1}) \; \geq \; b^{l} \; \geq \; length(\tilde{y}^{l}_{1}\tilde{\nu}^{m}_{1})\;;
\end{equation}
and also:

\begin{equation}
    \angle \tilde{\nu}^{l}_{1}\tilde{\tilde{b}}\tilde{\nu}^{k}_{1} \; > \; \angle \tilde{\nu}^{l}\tilde{\tilde{b}}\tilde{\nu}^{k} \; > \;  \angle \tilde{\nu}^{l}\tilde{\nu}^{m}\tilde{\nu}^{k} \; \geq \; \varphi_{1}\;;
\end{equation}
and so we dispose with cases 1) and 3) too, thus concluding the proof of part A).

\subsubsection{The extension to a fat $3$-dimensional triangulation} We start by observing that, if $u$ is a vertex of the simplex $s_{i}$ s.t. $f_{123} \cap s_i =
\tilde{\tau}$ -- where $f_{123}$ is a face of the tetrahedron $\sigma$ -- then the triangles $T^{kl}\,,\;T_{ij}$
and those produced by the subdivision of $Q_{kl}$ are also may also be conned from $u$.
\\We want to show that the simplices thus generated -- denoted by $V^{\delta}_{ij}\,=\,J(u_{\delta}\,,T_{ij})\,, \\ \delta\,=\,1,2,3,4$; etc. --
have big angles. We shall justify this affirmation for tetrahedra of type $V^{\delta}_{ij}$, the other cases being
completely analogous. Indeed:
\\ \qquad a) The angles $\tilde{\nu}^{k}_{i}\,,\:\tilde{\nu}^{k}_{j}$ and $\angle \tilde{\tilde{b}}$
are "big" by the very construction of $\triangle \tilde{\nu}^{k}_{i}\tilde{\nu}^{k}_{j}\tilde{\tilde{b}}$\,;
\\ \qquad b) The plane angles around the vertex $u_{\delta}$ are "big", since $length(\tilde{\nu}^{k}_{i}\tilde{\nu}^{k}_{j})\:\geq\:c_{1}b$\,,
and  because  $\tilde{\nu}^{k}_{i}\tilde{\nu}^{k}_{j}$  is included in the plane of $\tilde{\tau}$;
\\ \qquad c) The dihedral angles around $u_{\delta}$ are also "big" (by the argument above and by the "tetrahedral sinus formula").
To sum up, the angles denoted by $"L"$ in Fig. 13 are larger than some constant $\varphi_{2}$\,. However, the
basic type of small angles may still occur (see Fig. 14\,, where they are denoted by $"s"$).
\\In fact, the case of Fig. 3 (a) is associated with a low ratio "height/base side"; whereas, in this instance:

\begin{equation}
\frac{u_{\delta} \tilde{O}_{ij}}{\tilde{\nu}^{k}_{i}\tilde{\nu}^{k}_{j}}\;>\: \frac{u_{\delta}
\tilde{O}_{ij}}{\tilde{\nu}_{l} \tilde{\nu}_{m}}\;>\;c^{0}\;=\;const.\;;
\end{equation}
(by the fatness of the tetrahedron $u_{1}u_{2}u_{3}u_{4}$), so this case is excluded.
\\ Part of the angles covered by case (b) of Fig 14 are "big", for the following inequality (and its analogues hold)
\begin{equation}
\frac{u_{\delta}\tilde{\nu}^{k}_{i}}{\tilde{\nu}^{k}_{i}\tilde{\nu}^{k}_{j}}
\;>\;c^{\ast}\frac{u_{\delta}\tilde{\nu}_{k}}{\tilde{\nu}_{l}\tilde{\nu}_{k}}\:;\;\; c^{\ast}\:=\:const.\:;
\end{equation}

\begin{figure}[h]
\begin{center}
\includegraphics[scale=0.3]{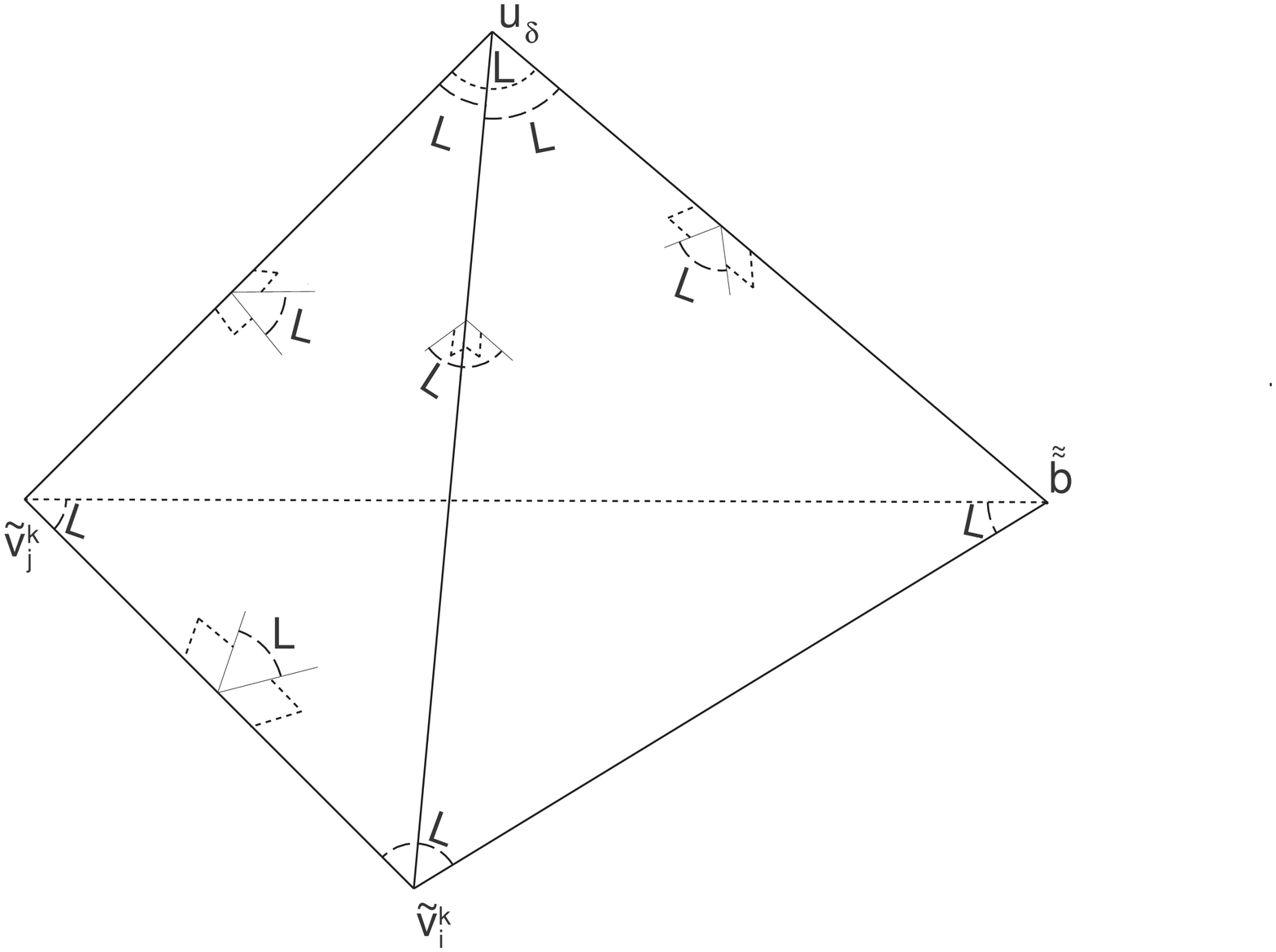}
\end{center}
\caption{}
\end{figure}

For a full proof of the "fatness" of $u_{\delta}\tilde{\nu}^{k}_{i}\tilde{\nu}^{k}_{j}\tilde{\nu}_{k}$ (to wit) we
still have to check the size of each of the following angles: $\angle
\tilde{\nu}^{k}_{i}\tilde{\nu}^{k}_{j}$\footnote{\,i.e. the dihedral angle between the faces $\triangle
u_{\delta}\tilde{\nu}^{k}_{i}\tilde{\nu}^{k}_{j}$ and
$\triangle\tilde{\nu}_{k}\tilde{\nu}^{k}_{i}\tilde{\nu}^{k}_{j}$.}, $\angle
\tilde{\nu}_{k}u_{\delta}\tilde{\nu}^{k}_{i}$ and $\angle \tilde{\nu}_{k}u_{\delta}\tilde{\nu}^{k}_{j}$\,.
\\To ensure the proper size of the last two angles one has to notice that the size of -- e.g.
$\angle \tilde{\nu}_{k}u_{\delta}\tilde{\nu}^{k}_{i}$ --  inversely proportional to:
\\ (a) the distance $\hat{\delta}$ from $u_{\delta}$ to $\tilde{\nu}^{k}_{l}\tilde{\nu}^{k}$\,,
\\and
\\ (b) the distance $\delta^{\ast}$ from $\bar{u}_{\delta}$ \footnote{\,Here $u_{\delta}\bar{u}_{\delta} \perp f_{123}$ } to $\tilde{\nu}^{k}_{l}\tilde{\nu}^{k}$\,.
\\ But, since the lengths $\tilde{\nu}^{k}_{l}\tilde{\nu}^{k}$ are bounded from below, there exists a minimal universal distance
$\delta^{\ast}_{0}$ s.t. if $\delta^{\ast} < \delta^{\ast}_{0}$ then  $\angle
\tilde{\nu}_{k}u_{\delta}\tilde{\nu}^{k}_{j} > \varphi^{\ast}_{2}$\,, for a suitable $\varphi^{\ast}_{2}$\,.

\begin{figure}[h]
\begin{center}
\includegraphics[scale=0.3]{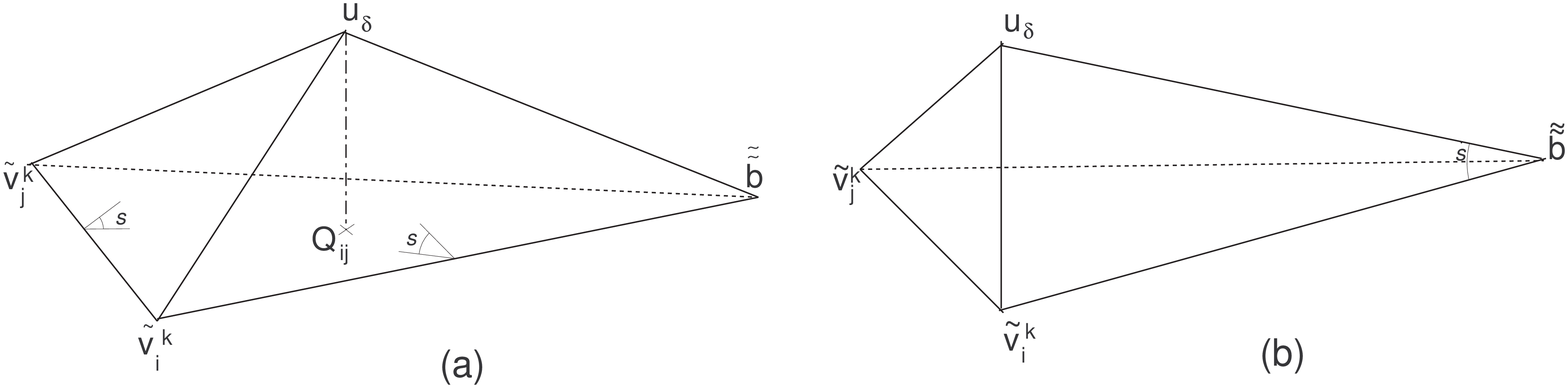}
\end{center}
\caption{}
\end{figure}

But, since the lengths $\tilde{\nu}^{k}_{l}\tilde{\nu}^{k}$ are bounded from below, there exists a minimal
universal distance $\delta^{\ast}_{0}$ s.t. if $\delta^{\ast} < \delta^{\ast}_{0}$ then  $\angle
\tilde{\nu}_{k}u_{\delta}\tilde{\nu}^{k}_{j} > \varphi^{\ast}_{2}$\,, for a suitable $\varphi^{\ast}_{2}$\,.
\\ By eventually decreasing $\delta^{\ast}_{0}$ to a new $\delta^{\ast}_{1}$\,, we can ensure that the angles of type $\angle
\tilde{\nu}_{k}u_{\delta}\tilde{\nu}^{k}_{i}$ are strictly grater than some $\varphi^{\ast \ast}_{2}$\,,
$\varphi^{\ast \ast}_{2} \leq \varphi^{\ast \ast}_{2}$\,, for each $\tilde{\nu}_{k} \in St(u_{\delta}) \cap
f_{123}$ and for each $\tilde{\nu}^{k}_{i}$, which concludes the proof of case(b). Indeed, we can vary the
position of $\bar{u}_{\delta}$ in the region $\Lambda_{\bar{u}_{\delta}}$\,, such that {\it all} the required
angles will be large -- see Fig. 15. (Such a small movement won't affect the fatness of the next stratum of
tetrahedra of type "$s$".)

\begin{figure}[h]
\begin{center}
\includegraphics[scale=0.3]{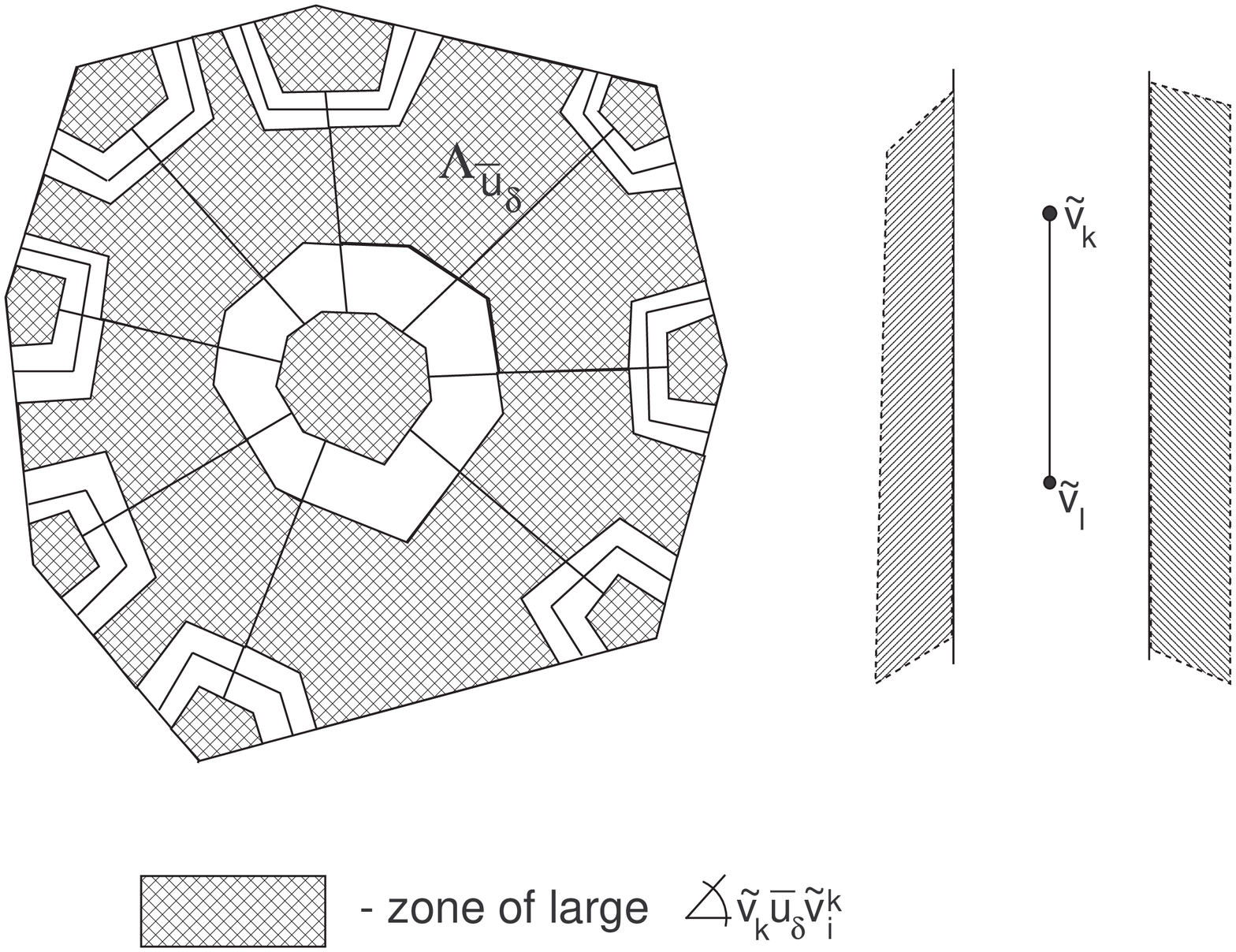}
\end{center}
\caption{}
\end{figure}

After disposing with this case, we can turn our attention to and to case (a)\,. Now, since $|\tan{\angle
\tilde{\nu}^{k}_{i}\tilde{\nu}^{k}_{j}}\,| = \frac{u_{\delta}\bar{u}_{\delta}}
{\bar{u}_{\delta}\bar{\bar{u}}_{\delta}}$\,, where  $u_{\delta}\bar{u}_{\delta} \perp
(\tilde{\nu}^{k}_{i}\nu_{k}\tilde{\nu}^{k}_{j})$, ${\bar{u_{\delta}}\bar{\bar{u}}_{\delta}} \perp
\tilde{\nu}^{k}_{i}\tilde{\nu}^{k}_{j}$\,, and since $\bar{u}_{\delta}\bar{\bar{u}}$ is known, we only have to
ensure that $u_{\delta}\bar{u}_{\delta}$ is bounded from below.

\begin{figure}[h]
\begin{center}
\includegraphics[scale=0.3]{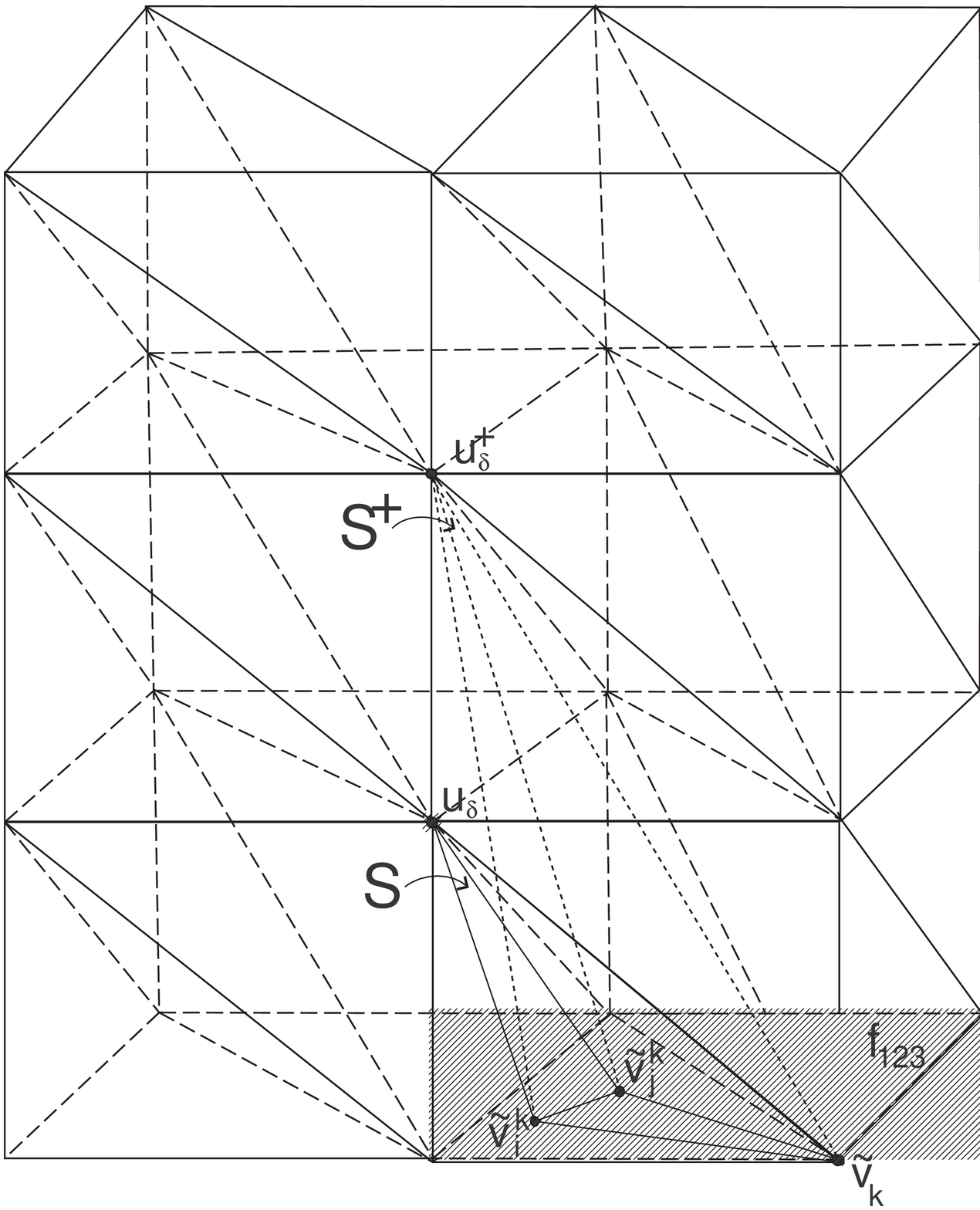}
\end{center}
\caption{}
\end{figure}

For this, one essentially makes use of the very specific form of the triangulation of the geometric neighbourhood
of an elliptic axes.
\\ Namely, instead of considering tetrahedra  $u_{\delta}\tilde{\nu}^{k}_{i}\tilde{\nu}^{k}_{j}\tilde{\nu}_{k}$\,, we
shall consider tetrahedra    $u^{+}_{\delta}\tilde{\nu}^{k}_{i}\tilde{\nu}^{k}_{j}\tilde{\nu}_{k}$\, where
$u^{+}_{\delta}$ is the vertex corresponding to $u_{\delta}$ in the next stratum of "$s$"-type tetrahedra, of the
family $\mathcal{S}$.\footnote{Of course -- as before -- we may have to subdivide the tetrahedra once or twice.}
\\Now the desired bound -- for $\frac{u^{+}_{\delta}\bar{u}^{+}_{\delta}}{\bar{u}^{+}_{\delta}\bar{\bar{u}}^{+}_{\delta}}$
instead of $\frac{u_{\delta}\bar{u}_{\delta}}{\bar{u}_{\delta}\bar{\bar{u}}_{\delta}}$ -- is achieved, whilst the
very existence of the constants obtained before isn't changed, only their magnitude; we shall denote them by
upper-right-superscript: e.g. $\varphi^{++}_{2}$.

\begin{rem}
The method we have just used is adaptable to the general context: consider -- instead of $u_{\delta}$ -- vertices
$u^{++}_{\delta}$\,, such that $dist_{eucl}(u_{\delta},u^{++}_{\delta}) = h_{0}$\,, where $h_{0} =
\frac{1}{3}h_{min}$\,, and $h_{min} =$ the minimal height of the simplices $s^{+} \in St(u_{\delta})$\,. The
points $u^{++}_{\delta}$ are to be chosen on the normal through $u_{\delta}$ to the plane $f_{123}$\,, so that the
combinatorics of the triangulation will suffer no alteration. Clearly, uniform bounds will again be attained.
\end{rem}

 In order to conclude the proof of Theorem 1.1 we still have to provide for a fat triangulation around the node
points. However, this missing case is easily dealt with by considering the intersections of the geometric
neighbourhoods of the elliptic axes. Such an intersection will automatically inherit from the tubular
neigbourhoods that generate it a {\it natural}, stratified, fat triangulation. (See Section \nolinebreak[4] 3.)
So, the arguments involved in the proof of the restricted case of un-intersecting axes do apply here, too. Thus we
conclude the
\[\hspace{6cm} \fbox{Proof of Theorem 1.1}\]

\section{Higher Dimensions}
 Our approach to the extension of our results to dimension $n \geq 4$ -- which we will expose only briefly -- is
based upon reduction to smaller dimensions. To do this, let us observe that if $S(n+1)$ is a simplex in $R^{n}$\,,
then any sectioning hyperplane separates the vertices of $S(n+1)$ into two groups of $p$ and $n+1-p$ vertices,
respectively. The notation will be as follows: $(p,n+1-p) \equiv (n+1-p,p)$ for the polytopes of the section; and
$(p\,|\,n+1-p)$ and $(n+1-p\,|\,p)$ for the resulting {\it frusta}.\footnote{\,See \cite{som} for notations and
precise definitions.}

The reduction to lower dimensions is permitted by:
\begin{lem}[\cite{som}]
    The frustum of type $(p\,|\,2)$  of $S(p+q)$ is isomorphic to the section of type $(p,q+1)$ of $S(p+q+1)$.
\end{lem}

The fattening algorithm is, in a nutshell, as follows: (i) Divide the polytopes obtained by sectioning into
simplices: first the section polytope, then the faces that are part of the original $S(n+1)$, while ensuring
fatness by the methods exposed in Section 3. If needed apply an inductive process on the dimension of the
simplices. (ii) "Fatten" the frusta by proving the existence of a locus of points where from {\it all} the
$(n-k)$-faces are seen at big $(n-k)$-dimensional angles.\footnote{\,See \cite{som}  for the relation between the
$p$-dimensional angles of a $n$-simplex, $p = 2,\ldots,n-1$.}
\\ Theorem 1.1 now follows. However, we have to understand much better the geometry of the elliptic locus of a Kleinian group with torsion, acting upon
$H^{n},\; n \geq 4$\,. We are, however, fortunate, for the fixed set of an elliptic transformation is a
$k$-dimensional hyperbolic plane, $0 \leq k \leq n-2$\,; thus providing the fixed point set with a geometric
neighbourhood -- together with its natural fat triangulation. Some complications may arise because different
elliptics may well have fixed loci of different dimensions\footnote{\,For a detailed discussion on the diversity
of the elliptics (and their fixed loci) see \cite{ap}.}, so the respective simplices will have different
dimensions and fatnesses. However this is easy to remedy by completing them to $n$-dimensional simplices, by
"expanding" the low dimensional neighbouhoods to maximal dimension in a product manner.

\section{Appendix -- Sketch of Proof of Theorem 2.4.} Since in dimension $3$ parabolic elements may have only rank $1$ or $2$ (See \cite{ms}, \cite{abi1}),
suffice to analyze only the following two cases:
\\ (A) If $G$  contains no parabolic elements or if $G$ contains only rank $1$
parabolics (see\cite{p}) then, by a theorem of Milnor-Gromov (see \cite{gro}, \cite{gdlh}), $\pi(G)$ is
word-hyperbolic and it follows, by a corollary of Rips' theorem (see\cite{gdlh}), that it contains only a finite
number of classes of elements of finite order, i.e. elliptics.
\\ (B) If $G$ contains parabolics of rank $2$, then one one can opt for one of the
following paths:
\\ (a) Use Scott's Theorem (See \cite{abi1}) that states that $3$-manifolds are {\it compact core manifolds}\footnote{\,i.e. if $M$ is a $3$-manifold s.t. $\pi_{1}(M)$ is finitely generated, then there exists $M_{0} \subset
M$, $M_{0}$ compact and s.t. $\pi_{1}(M_{0}) \cong \pi_{1}(M_{1})$}, then apply the Milnor-Gromov Theorem for
$M_{0}$.
\\ or
\\ (b)"Double" the manifold with respect to its cuspidal ends (see \cite{mor},
\cite{abi1}\,).
\\Since the number of cusps is finite (by a theorem of Sullivan's  -- see \cite{abi1}, \cite{kp})
and since the double $\widetilde{M}$ is compact it follows, exactly as above, that $\pi(\widetilde{M})$ has only a
finite number of conjugacy classes of elliptics, hence so has $\pi_{1}(M)$.
\\\hspace*{\fill}$\Box$ \par

\begin{thebibliography}{99}
\bibitem[Abi1]{abi1}
Abikoff, W.\,: {\it Kleinian groups -- geometrically finite and geometrically perverse }, Contemporary
Mathematics, Vol. 74, 1988, pp. 1-50.

\bibitem[Abi2]{abi2}
Abikoff, W.\,: {\it Kleinian Groups}, Lecture Notes, given at "The Technion -- Israel Institute of Technology",
Haifa, Israel, 1996-1997.

\bibitem[AH]{ah}
Abikoff, W. and  Hass, A.H.\,: {\it Nondiscrete groups of hyperbolic motions}, Bull. London Math. Soc. 22, 1990,
no.3, pp. 233-238.

\bibitem[Al]{al}
Alexander, J.W.\,: {\it Note on Riemmann spaces}, Bull. Amer. Math. Soc. 26, 1920, pp. 370-372.

\bibitem[AP]{ap}
Apanasov, B.N. \,: {\it Klein Groups in Space}, Sib. Math. J. 16, 1975, pp. 679-684.

\bibitem[BAC]{bac}
Br\^{a}nzei, D.\,, Ani\c{t}a, S. and Cocea, C.\,: {\it Planul si spa\c{t}iul euclidian}, Editura Academiei R.S.R.,
Bucure\c{s}ti, 1986.

\bibitem[Bea]{bea}
Beardon, A. F.\,: {\it The Geometry of Discrete Groups}, Springer Verlag, GTM 91, N.Y., 1982.

\bibitem[BM]{bm}
Bowditch, B.H. and Mess, G.\,: {\it A 4-Dimensional Kleinian Group}, Transaction of the Amer. Math. Soc., Vol.
344, No.1, 1994, pp. 390-405.

\bibitem[BrM]{brm}
Brooks, R. and Matelski, J.P.\,: {\it Collars in Kleinian groups}, Duke Math.J.\,, Vol. 49, No. 1, 1982, pp.
163-182.


\bibitem[Ca1]{ca1}
Cairns, S.S.\,: {\it On the triangulation of regular loci}, Ann. of Math. 35, 1934,
 pp. 579-587.

\bibitem[Ca2]{ca2}
Cairns, S.S.\,: {\it Polyhedral approximation to regular loci}, Ann. of Math. 37, 1936, pp. 409-419.

\bibitem[Ca3]{ca3}
Cairns, S.S.\,: {\it A simple triangulation method for smooth manifolds}, Bull. Amer. Math. Soc. 67, 1961, pp.
380-390.

\bibitem[CEG]{ceg}
Canary, R.D.\,, Epstein, D.B.A. and Green, P. \,: {\it Notes on the notes of Thurston}, London Math. Soc. Lecture
Notes Series 111, 1987, pp. 3-92.

\bibitem[Cox]{cox}
Coxeter, H. S. M. \,: {\it Regular Polytopes}, Second Edition, Macmillan, NY, 1963.
\bibitem[DM]{dm}
Derevin, D.A. and Mednikh, A.D.\,: {\it Geometric proprieties of Discrete groups acting with fixed points in
Lobachevsky space}, Soviet Math. Dokl., Vol. 37, 1988, No.3, pp. 614-617.

\bibitem[Ep$^{++}$]{ep}
Epstein, D.B.A.\, et.al (editor)\,: {\it Word Processing in groups},
 Jones and Bartlett, Boston, MA, 1992.
\bibitem[FM]{fm}
Feighn, M. and Mess, G.\,: {\it Conjugacy classes of finite subgroups of Kleinian groups}, Amer. J. of Math., 113,
1991, pp. 179-188.
\bibitem [Fe]{fe}
Ferraroti, M.\,: {\it Triangulatione analitica per variet\`{a}}, Bolletino U.M.I. (5) {\bf 17}-A (1980), 79-84.

\bibitem[GdlH]{gdlh}
Guys, E. and  de la Harpe, P. (ed.) \,: {\it  Sur les Groupes Hyperbolique apres Mikhael Gromov}, Progr. Math.,
Vol 83., Birkh\"{a}user, Boston, MA, 1990.

\bibitem[GM1]{gm1}
Gehring, F.W.,  Martin G.J.\,, {\it Commutators, collars and the geometry of M\"{o}bius groups}, J. Anal. Math.
Vol. 63, 1994, pp. 174-219.

\bibitem[GM1]{gm2}
Gehring, F.W., Martin, G.J.\,, {\it On the Margulis constant for Kleinian groups, I}, Ann. Acad. Sci. Fenn., Vol.
21, 1996, pp. 439-462.

\bibitem[GMMR]{gmmr}
Gehring, F.W., Maclachlan, C., Martin, G.J., and Reed, A.W.\,: {\it Aritmecity, Discreteness and Volume}, Trans.
Amer. Math. Soc. 349, 1997, 3611-3643.

\bibitem[Gro]{gro}
Gromov, M.\,: {\it Infinite groups as geometric objects}, in "Geometric Group Theory", (Niblo, G.A. and Roller,
M.A. ed.), Spriger Verlag, MSRI Publ. 8, 1987, pp. 75-263.


\bibitem[H]{h}
Hamilton, E.\,: {\it Geometrical finiteness for hyperbolic orbifolds}, Topology, Vol.37, No.3, 1998,pp. 635-637.


\bibitem[Hu]{hu}
Hudson, J.F. \,: {\it Piecewise Linear Topology}, Math. Lect. Notes Series, Benjamin, N.Y., 1969.

\bibitem[J]{j}
J\o rgensen, T.\,: {\it On discrete groups of M\"{o}bius transformations}, Amer. Journ. of Math., Vol. 98, No.3,
pp.739-749.


\bibitem[KP]{kp}
Kapovitch, M.E. and Potyagailo, L.\,: {\it On the Absence of Ahlfors and Sullivan  theorems for Kleinian groups in
higher dimensions}, Sib. Math. J., Vol. 32, No. 1, 1991, pp. 227-237.




\bibitem[Med]{med}
Mednikh, A.D.\,: {\it Automorphism groups of the three-dimensional hyperbolic manifolds}, Soviet Math. Dokl., Vol.
32, 1985, No.3, pp. 633-636.

\bibitem[Mor]{mor}
Morgan, J.W., \,: {\it On Thurston's Uniformization Theorem for Three-Dimensional Manifolds}, in "The Smith
Conjecture", (Morgan, J.W. and Bass, H. ed.), Academic Press, N.Y., 1984, pp. 37-126.

\bibitem[MS1]{ms1}
Martio, O., and Srebro, U. \,: {\it Automorphic quasimeromorphic mappings in $R^{n}$} , Acta Math. 195, 1975, pp.
221-247.

\bibitem[MS2]{ms2}
Martio, O., and  Srebro, U.\,: {\it On the existence of automorphic quasimeromorphic mappings in $R^{n}$}, Ann.
Acad. Sci. Fenn., Series\,I Math., Vol. 3, 1977, pp. 123-130.

\bibitem[Ms]{ms}
Maskit, B.\,: {\it Kleinian Groups}, Springer Verlag, GDM 287, N.Y., 1987.

\bibitem[Mun]{mun}
Munkres, J. R.\,: {\it Elementary Differential Topology}, (rev. ed.) Princeton University Press, Princeton, N.J.,
1966.

\bibitem[NW]{nw}
Nicholls, P.J. and Waterman, P.L. \,: {\it The boundary of convex fundamental domains for Fuchsian groups}, Ann.
Acad. Sci. Fenn., Ser A I Math. 15, !990, no. 1, pp. 1-25.


\bibitem[P]{p}
Potyagailo, L. \,: {\it Finitely generated Kleinian groups in $3$-space and $3$-manifolds of infinite homotopy
type}, Trans. of Amer. Math. Soc., Vol. 344, No. 1, 1994, pp. 57-77.


\bibitem[Pe]{pe}
Peltonen, K.\,: {\it On the existence of quasiregular mappings}, Ann. Acad. Sci. Fenn., Series\,I Math.,
 Dissertationes, 1992.

\bibitem[Rat]{rat}
Ratcliffe, J.G. \,: {\it Foundations of Hyperbolic Manifolds}, GTM 194, Springer Verlag, N.Y., 1994.

\bibitem[S1]{s1}
Saucan, E. \: {\it On the existence of quasimeromorphic mappings}, in preparation.

\bibitem[S2]{s2}
Saucan, E. \: {\it A note on a theorem of Munkres}, in preparation.

\bibitem[S3]{teza}
Saucan, E. \: {\it Ph. Thesis, technion, in preparation}


\bibitem[Som]{som}
Sommerville, D.M.Y.\,: {\it An introduction to the Geometry of  $N$ Dimensions}, Dover Publications, N.Y., 1958.



\bibitem[Spi V]{spi5}
Spivak, M.,\,:  {\it A comprehensive Introduction to Differential Geometry, volume V}, Publish or Perish, Boston,
MA, 1975.

\bibitem[Sr]{sr}
Srebro, U. \,: {\it Non-existence of  Automorphic Quasimeromorphic Mappings}, Analysis and topology, World Sci.
Publishing, River Edge, Nj, 1998.

\bibitem[Th1]{th1}
Thurston, W.\,: {\it The Geometry and Topology of $3$-Manifolds}, preliminary version, chap. 1-9, 1990.

\bibitem[Th2]{th2}
Thurston, W. \,: {\it Three-Dimensional Geometry and Topology, vol.1, (Edited by S. Levy)}, Princeton University
Press, Princeton, N.J. 1997.

\bibitem[Tu]{tu}
Tukia, P.\,: {\it Automorphic Quasimeromorphic Mappings for Torsionless Hyperbolic Groups}, Ann. Acad. Sci. Fenn.,
10, 1985, pp. 545-560.



\bibitem[Wh]{wh}
Whitehead, J.H.C.\,: {\it On $\mathcal{C}^{1}$-complexes}, Ann. of Math., 41 (1940), 809-824.

\end{thebibliography}

\end{document}